\documentclass[11pt,a4paper]{amsart}
\usepackage[dvips]{epsfig}
\usepackage{graphics}
\usepackage{latexsym}
\usepackage{verbatim}
\usepackage{amsmath}
\usepackage{amsthm}
\usepackage{amssymb}
\usepackage{hyperref}
\newtheorem{theorem}{Theorem}[section]

\newtheorem{proposition}[theorem]{Proposition}

\newtheorem{remark}[theorem]{Remark}

\newcommand{\M}{\mathcal{M}}
\newcommand{\HH}{\mathcal{H}}

\newcommand{\R}{\mathbb{R}}

\newcommand{\Tau}{\mathcal{T}}

\newcommand{\QQ}{\mathbb{Q}}

\newcommand{\Q}{\mathcal{Q}}
\newcommand{\G}{\mathcal{G}}

\def\var{\varepsilon}

\def\signff{\bigskip\bigskip\hspace{80mm}
\vbox{{\sc Francis Filbet\par\vspace{3mm}
Universit\'e de Lyon,\par
Universit\'e Lyon I, CNRS \par
UMR 5208, Institut Camille Jordan \par
43, Boulevard du 11 Novembre 1918\par
69622 Villeurbanne cedex, FRANCE\par\vspace{3mm}
e-mail:} filbet@math.univ-lyon1.fr }}

\def\signsj{\bigskip\bigskip\hspace{80mm}
\vbox{{\sc Shi Jin\par\vspace{3mm}
 Department of Mathematics \par
University of Wisconsin\par
Madison, WI 53706, USA 
\par\vspace{3mm} e-mail:} jin@math.wisc.edu }}

\title[An Asymptotic Preserving Scheme for the ES-BGK model]{An Asymptotic Preserving Scheme for the ES-BGK model of the Boltzmann equation}\thanks{F. Filbet is partially supported by the European Research Council ERC Starting Grant 2009,  project 239983-\textit{NuSiKiMo}. S. Jin was partially supported by NSF grant No. DMS-0608720, NSF FRG grant DMS-0757285, and a Van Vleck Distinguished Research Prize from University of Wisconsin-Madison.}

\author{Francis Filbet \and Shi Jin }

\begin{document}

\maketitle

\begin{abstract}
In this paper, we study a time discrete scheme for the initial value problem 
of the ES-BGK  kinetic equation. Numerically solving these equations are 
challenging due to the nonlinear stiff collision (source) terms induced by 
small mean free or relaxation time. We study an implicit-explicit (IMEX)
time discretization in which the convection is explicit while the relaxation
term is implicit to overcome the stiffness.
We first show how the implicit relaxation can be solved explicitly, and then
prove asymptotically that this time discretization drives the density 
distribution toward the local Maxwellian when the mean free time goes
to zero while the numerical time step is held fixed. This  naturally imposes 
an asymptotic-preserving scheme in the Euler limit. The scheme so designed  
does not need any nonlinear iterative solver for the
implicit relaxation term. Moreover, it can capture the macroscopic fluid 
dynamic (Euler) limit even if the small scale determined by the  Knudsen 
number is not numerically resolved. We also show that it is consistent to the 
compressible Navier-Stokes equations if the viscosity and heat conductivity are numerically resolved.  Several numerical examples, in both one and two space
dimensions, are used to demonstrate the desired behavior of this scheme.
\end{abstract}

\tableofcontents

\section{Introduction}\label{sec1}
\setcounter{equation}{0}
When gas is in thermal non-equilibrium, which are encountered frequently in 
hypersonic flows, vehicles at high altitudes and flows expanding into vacuum, 
the macroscopic constitutive laws based on the continuum hypothesis tend to 
breakdown.  A critical parameter that characterizes the rarifiedness of 
the gas is the Knudsen number ($\varepsilon = \lambda/L$), where $\lambda$ is 
the average distance traveled by the molecules between collisions, or the mean free path, and $L$ is the characteristic length scale. 
When the flow gradients are large, such as in shock or boundary layers,
continuum fluid dynamics equations are not adequate, and one needs to use
a kinetic equation. 
The fundamental kinetic equation for rarefied gas is  the Boltzmann equation,
\begin{equation}
\label{eq:1}
\frac{\partial f}{\partial t} + v\,\cdot\,\nabla_x f \,=\, \frac{1}{\varepsilon}\,\Q(f),	
\end{equation}
which governs the evolution of the density $f(t, x, v)$ of monoatomic 
particles 
in the phase, where $(x,v)\in\R^{d_x}\times \R^{d_v}$. Boltzmann's collision 
operator has the fundamental properties of conserving mass, momentum and energy: at the formal level
 \begin{equation*}
 \int_{{\R}^d}Q(f,f) \, \phi(v)\,dv = 0, \qquad
 \phi(v)=1,v,|v|^2,
 \end{equation*}
Moreover, the equilibrium is the local Maxwellian distribution
 \begin{equation}
\label{maxw}
 \M[f](v)=\frac{\rho}{(2\pi T)^{d_v/2}}
 \exp \left( - \frac{\vert u - v \vert^2} {2T} \right), 
 \end{equation}
where $\rho,\,u,\,T$ are the {\em density}, {\em macroscopic velocity}
and {\em temperature} of the gas, defined by
\begin{eqnarray}
\label{ru}
\rho = \int_{{\R}^{d_v}}f(v)\,dv=\int_{{\R}^{d_v}}\M[f](v), \quad u =
 \frac{1}{\rho}\int_{{\R}^{d_v}}v\,f(v)\,dv=
 \frac{1}{\rho}\int_{{\R}^{d_v}}v\,\M[f](v)\,dv,
\end{eqnarray}
and
\begin{eqnarray}
\label{t} 
T = \frac{1}{d_v\rho}
 \int_{{\R}^{d_v}}\vert u - v \vert^2\,f(v)\,dv=\frac{1}{d_v\rho}
 \int_{{\R}^{d_v}}\vert u - v \vert^2\,\M[f](v)\,dv.
 \end{eqnarray}

The Boltzmann equation is closely related to the Navier-Stokes system which governs the evolution of macroscopic density, momentum and energy in the
continuum regime:
\begin{equation}
\label{eq:CNS}
\left\{
\begin{array}{l}
\displaystyle \frac{\partial \rho}{\partial t} \,+\,{\rm div}_x \left(\rho \,u\right) \,=\, 0,
\\
\,
\\
\displaystyle \frac{\partial \rho \,u }{\partial t} \,+\,{\rm div}_x \left(\rho \,u \otimes u \,+\, p \,{\rm I}\right) \,=\, \varepsilon\,{\rm div}_x[\mu\,\sigma(u)],
\\
\,
\\
\displaystyle \frac{\partial E }{\partial t} \,+\,{\rm div}_x \left( (E\,+\, p)\,u\right)\,=\, \varepsilon\, {\rm div}_x\left(\mu\,\sigma(u)\,u \,+\, \kappa\,\nabla_x T\right).
\end{array}\right.
\end{equation}
where $p$ is the pressure, $E$ represents the total energy
$$
E \,\,=\,\, \frac{1}{2} \,\rho \,u^2 \,+\, \frac{d_v}{2}\,\rho\, T,
$$
and ${\rm I}$ is the identity matrix. Moreover, the tensor $\sigma(u)$ denotes the strain rate tensor given by
$$
\sigma(u) = \left( \nabla_x u + (\nabla_x u)^t \right) \,-\,\frac{2}{d_v} {\rm div}_x u \,{\rm I}.
$$
These equations  constitute a system of $2 + d_v$ equations in $3 + d_v$ unknowns. The pressure  is related to the internal energy by the constitutive relation for a polytropic gas  
$$
p \,=\, (\gamma - 1)\,\left( E \,-\,\frac{1}{2}\,\rho \,u^2\right),
$$ 
where the polytropic constant $\gamma = (d_v + 2)/d_v$ represents the ratio between specific heat at constant pressure and at constant volume, thus yielding
 $p=\rho\,T$ 
while the viscosity $\mu =\mu(T)$ and the thermal conductivity $\kappa =\kappa(T)$ are defined according to the linearized Boltzmann operator with respect to the local Maxwellian \cite{bgl:91}.

Since the quadratic collision operator $\Q(f)$ has a rather complex form, 
simpler models have been introduced. The main requirement is to build models 
that have the right conservations, entropy condition described by the
$H$-theorem, and have the fluid dynamics (Euler and Navier-Stokes)
limits with the correct transport coefficients.  The simplest model is the 
so-called BGK model introduced by Bhatnagar, Gross and Krook \cite{BGK}. It is based on relaxation towards the local Maxwellian
\begin{equation}
\label{C-BGK}
\Q(f)= \frac{\tau}{\varepsilon}\,\left(\M[f] \,-\, f\right),
\end{equation}
where $\tau$ depends on  macroscopic quantities $\rho$, $u$ and $T$.
 
This model conserves mass, momentum and total energy, and has the 
correct Euler limit when $\varepsilon \to 0$. But in the Chapman-Enskog 
expansion, the transport coefficients, that is $\mu$ and $\kappa$ obtained at the Navier-Stokes level are not satisfactory. In particular, the Prandtl number defined by
$$
{\rm Pr }\,=\, \frac{\gamma}{\gamma-1} \,\frac{\mu}{\kappa}\,,
$$
which relates the viscosity to the heat conductivity, is equal to $1$, whereas for most gases, we have Pr $< 1$.  For instance, the hard-sphere model for a monoatomic gas ($\gamma = 5/3$) in the Boltzmann equation leads to a Prandtl number very close to $2/3$. 


One model, proposed by Holway \cite{cite11}, has the desired property of
having the correct conservation laws, yields the Navier-Stokes
approximation via the Chapman-Enskog expansion with a Prandtl number less 
than one, and yet is endowed with the
entropy condition  \cite{Perth}. See also \cite{cite6, Cer}. This model
is referred to as the ellipsoidal statistical model (ES-BGK),
where the Maxwellian $\M[f]$ in the relaxation term of (\ref{C-BGK})
is replaced by an anisotropic Gaussian $\G[f]$. In order to introduce the 
Gaussian model, we need further notations. Define the opposite of the stress tensor
\begin{equation}
\label{thet}
\Theta(t,x)\, =\, \frac{1}{\rho} \int_{{\R}^{d_v}}  (v-u)\otimes (v-u)\,f(t,x,v) \,dv.
\end{equation} 
Therefore the translational temperature is related to the $T = {\rm tr}(\Theta)/d_v$. We finally introduce the corrected tensor
\begin{equation}
\label{tau}
\Tau(t,x) \,\,=\,\, \left[(1-\nu) \, T \,{\rm I} \,\,+\,\,\nu \,\Theta\right](t,x),  
\end{equation} 
which can be viewed as a linear combination of the initial stress tensor $\Theta$ and of the isotropic stress tensor $T \,{\rm I}$ developed by a Maxwellian distribution.

The Gaussian model introduces a corrected BGK collision operator by replacing the local equilibrium Maxwellian by the Gaussian $\G[f]$ defined by
$$
\G[f]= \frac{\rho}{\sqrt{{\rm det}(2\pi\,\Tau)}}\,\exp\left(-\frac{(v-u)\,\Tau^{-1}\,(v-u)}{2}\right).
$$ 
Thus, the corresponding collision operator is now
\begin{equation}
\Q(f)= \frac{\tau}{\varepsilon} \left( \G[f] - f\right),
\label{eq:esbgk}
\end{equation}
where $\tau$ depends on $\rho$, $u$ and $T$,  the parameter $-1/2 \leq \nu < 1$ is used to modify the value of the Prandtl number through the formula \cite{Perth}
$$
\frac{2}{3} \,\leq\, {\rm Pr}  \,=\, \frac{1}{1-\nu} \,\leq\, +\infty\quad {\rm for } \quad\nu\in [-1/2\,,\,1].
$$
It first follows from the above definitions that 
$$
\left\{\begin{array}{ll}
\displaystyle{\int_{\R^{d_v}} f(v)\,dv  = \int_{\R^{d_v}} \G[f](v)\,dv  = \rho,} &
\\
\;
\\
\displaystyle{ \int_{\R^{d_v}}v\, f(v)\,dv = \int_{\R^{d_v}} v\,\G[f](v)\,dv  = \rho\,u,} &
\\
\,
\\
\displaystyle{ \int_{\R^{d_v}} \frac{|v|^2}{2}\,f(v)\,dv  = \int_{\R^{d_v}} \frac{|v|^2}{2}\,\G[f](v)\,dv  = E} &
\end{array}\right.
$$
and
$$
\left\{
\begin{array}{ll}
\displaystyle \int_{\R^{d_v}} (v-u) \otimes (v-u)\, f(v)\,dv \,=\, \rho\, \Theta,&
\\
\,
\\
\displaystyle\int_{\R^{d_v}} (v-u) \otimes (v-u)\, \G[f] \,dv \,
=\, \rho \,\Tau.&
\end{array}\right.
$$
This implies that this collision operator does indeed conserve mass, momentum and energy as imposed. The collision frequency, $\nu$, involves the Prandtl 
number Pr as a free parameter. This allows the ES-BGK collision model to reproduce transport coefficients, viscosity and thermal conductivity, in the Chapman-Enslog expansion \cite{bgl:91,bouchut:golse}, recovering
the Navier-Stokes equations density $\rho$, momentum $\rho\,u$ and temperature $T$, with the correct Prandtl number.

In this paper we study a temporally implicit-explicit (IMEX) discretization
of the ES-BGK model. The advantage of such a time discretization is that
it is uniformly stable with respect to the small Knudsen number, thus
removing the stiffness of the relaxation term, yet the implicit
relaxation term can be solved explicitly, thanks to the special structure
of the relaxation term. Although such a property was realized for the
classical BGK operator \cite{CP}, the ES-BGK operator is different, and
we realized that one has to compute to the higher moment in order to evaluate
the implicit Gaussian distribution explicitly.  Furthermore, we
show that this time discretization is asymptotic-preserving \cite{J},
an important property for the scheme to be robust in the fluid dynamic
regime, allowing it to capture the fluid dynamic behavior without
resolving numerically the small Knudsen number. We further show that this
discretization is consistent to the compressible Navier-Stokes equations
with the correct Prandtl number, if the viscosity and heat flux terms are
suitably resolved numerically. We then validate this numerical method
by presenting
numerical results based on this scheme, and compare them with the solutions
of the Boltzmann equation and the corresponding Navier-Stokes equations,
in both one and two space dimensions.

\section{An Asymptotic Preserving scheme to the ES-BGK equation}
\label{sec3}
\setcounter{equation}{0}

Past progress on developing robust numerical schemes for  kinetic equations 
that also work in the fluid regimes has been guided by the fluid dynamic limit,
in the framework of {\it asymptotic-preserving} (AP) scheme. As summarized
by Jin
\cite{J}, a scheme for the kinetic equation is AP if

\begin{itemize}
\item  it preserves the
discrete analogy of the Chapman-Enskog expansion, namely, it is a
suitable scheme for the kinetic equation, yet, when holding
the mesh size and time step fixed and letting
the Knudsen number go to zero, the scheme becomes a suitable scheme
for the limiting Euler equations
\item implicit collision terms can be implemented explicitly, or at least
more efficiently than using the Newton type solvers for nonlinear
algebraic systems.
\end{itemize}

We now introduce the time discretization for the ES-BGK equation (\ref{eq:1}), (\ref{eq:esbgk})
\begin{equation} 
\label{eq:deb}
\left\{
\begin{array}{l}
  \displaystyle{\frac{\partial f}{\partial t}  \,+\, v\,\nabla_x f \,=\,  
\frac{\tau}{\varepsilon}\,(\G[f] -f),  
   \quad x \in \Omega\subset\R^{d_x},\, v\in \R^{d_v},} 
  \\
  \,
  \\
  f(0,x,v)  \,=\, f_{0}(x,v), \quad x\in\Omega,\,v\in\R^{d_v}, 
\end{array}\right.
\end{equation}
where $\tau$ depends on $\rho$, $u$ and $T$. 

The time discretization is an implicit-explicit (IMEX) scheme.
 Since the convection term in (\ref{eq:deb}) is not stiff, we will
treat it explicitly. The source terms on the right hand side of 
(\ref{eq:deb}) will be handled using an implicit solver. We simply apply a first order implicit-explicit (IMEX) scheme, 
\begin{equation} 
\label{sch:perturb}
\left\{
\begin{array}{l}
  \displaystyle{\frac{f^{n+1}-f^n}{\Delta t }  + v\cdot\nabla_x f^n \,=\, 
\frac{\tau^{n+1}}{\varepsilon}(\G[f^{n+1}]-f^{n+1}),}
  \\
  \,
  \\
  f^0(x,v)  = f_{0}(x,v)\,.
\end{array}\right.
\end{equation}
This can be written as
\begin{eqnarray}
\label{AP-1}
f^{n+1} &=& \frac{\varepsilon}{\varepsilon+\tau^{n+1}\Delta t} \left[f^n - \Delta t  \,v\cdot\nabla_x f^n\right] \,+\, \frac{\tau^{n+1} \Delta t}{\varepsilon+\tau^{n+1} \Delta t}\,\G[f^{n+1}],
 \end{eqnarray}
where $\G(f^{n+1})$ is the anisotropic Maxwellian distribution computed from $f^{n+1}$. Although (\ref{AP-1}) appears  nonlinearly implicit, since the computation of $f^{n+1}$ requires the knowledge of $\G[f^{n+1}]$, it can be solved explicitly. Specifically, upon multiplying (\ref{AP-1}) by $\phi(v)$ 
defined by
$$
\phi(v) \,:=\, \left(1,v,\frac{|v|^2}{2}\right) 
$$
and use the conservation properties of $\Q$ and the definition of $\G[f]$ in  (\ref{ru}), (\ref{t}), we define the macroscopic quantity $U$ by  $U:=(\rho,\rho\,u,E)$ computed from $f$ and get \cite{CP,PP07}
$$
U^{n+1} = \frac{\varepsilon}{\varepsilon+\tau^{n+1} \Delta t} \int_{\R^{d_v}} \phi(v)\,(f^n-\Delta t\,  v \cdot \nabla_x f^n) \, dv 
+ \frac{\tau^{n+1}\Delta t}{\varepsilon+\tau^{n+1} \Delta t} \,\int_{\R^{d_v}} \phi(v)\G[f^{n+1}](v)\,dv,
$$
or simply
\begin{equation}
\label{AP-11}
U^{n+1} =  \int_{\R^{d_v}}  \phi(v)\,(f^n-\Delta t v \cdot \nabla_x f^n)  \, dv \,.
\end{equation}
Thus $U^{n+1}$ can be obtained explicitly. This gives $\rho^{n+1},
u^{n+1}$ and $T^{n+1}$.
Unfortunately, it is not enough to define $\G[f^{n+1}]$ for which we need $\rho^{n+1}\,\Theta^{n+1}$. Therefore, we define the tensor $\Sigma$ by
\begin{equation}
\label{sigm}
\Sigma^{n+1} := \int_{\R^{d_v}} v \otimes v\, f^{n+1}\,dv \,= \;\rho^{n+1}\, \left(\Theta^{n+1} \,+\, u^{n+1} \otimes u^{n+1}\right)
\end{equation}
and multiply the scheme (\ref{AP-1}) by $v \otimes v$. Using the fact that 
$$
\int_{\R^{d_v}} v \otimes v\, \G[f](v)\,dv =  \rho\, \left(\Tau\,+\, u \otimes u\right),
$$
and (\ref{sigm}), we get that
\begin{eqnarray}
\label{AP-111}
\Sigma^{n+1} &=& \frac{\varepsilon}{\varepsilon+(1-\nu)\,\tau^{n+1}\,\Delta t}\left(\Sigma^n \,-\, \Delta t\,\int_{\R^{d_v}} v \otimes v \,v \cdot \nabla_x f^n dv \right)
\\
\,
\nonumber
\\
&& +\, \frac{(1-\nu)\,\tau^{n+1}\,\Delta t}{\varepsilon+(1-\nu)\,\tau^{n+1}\,\Delta t} \, \rho^{n+1}\,\left(T^{n+1} \,{\rm Id} \,+\, u^{n+1}\otimes u^{n+1} \right)   
\nonumber
\end{eqnarray}
Now $\G[f^{n+1}]$ can be obtained explicitly from $U^{n+1}$ and $\Sigma^{n+1}$ and then $f^{n+1}$ from  (\ref{AP-1}).

Finally the scheme reads
\begin{equation}
\label{AP-final}
\left\{
\begin{array}{lll}
U^{n+1} &=&  \displaystyle{\int_{\R^{d_v}}  \phi(v)\,(f^n-\Delta t v \cdot \nabla_x f^n)  \, dv,} 
\\
\,
\\
\Sigma^{n+1} &=& \displaystyle{\frac{\varepsilon}{\varepsilon+(1-\nu)\,\tau^{n+1}\,\Delta t}\left(\Sigma^n \,-\, \Delta t\,\int_{\R^{d_v}} v \otimes v \,v \cdot \nabla_x f^n dv \right)}
\\
\,
\\
&&\displaystyle{ +\, \frac{(1-\nu)\,\tau^{n+1}\,\Delta t}{\varepsilon+(1-\nu)\,\tau^{n+1}\,\Delta t} \, \rho^{n+1}\,\left(T^{n+1} \,{\rm Id} \,+\, u^{n+1}\otimes u^{n+1} \right).}
\\
\,
\\
f^{n+1} &=& \displaystyle{\frac{\varepsilon}{\varepsilon+\tau^{n+1}\Delta t} \left[f^n - \Delta t  \,v\cdot\nabla_x f^n\right] \,+\, \frac{\tau^{n+1} \Delta t}{\varepsilon+\tau^{n+1} \Delta t}\,\G[f^{n+1}]},   
\end{array}\right.
\end{equation}

In summary, although (\ref{sch:perturb}) is nonlinearly
implicit, it can be solved {\it explicitly}, thus satisfies the second condition of an AP scheme.

Now let us prove that the scheme (\ref{AP-final}) preserves the correct asymptotic.

\begin{remark}
When the IMEX scheme (\ref{sch:perturb}) is applied to the classical
BGK equation (\ref{C-BGK}), $U^{n+1}$ in (\ref{AP-11}) will completely
defines $\M[f^{n+1}]$ given by (\ref{maxw}) \cite{CP, FJ}). The steps after equation
(\ref{AP-11}) are new ideas introduced in this paper for the ES-BGK
model.
\end{remark}
    
\begin{proposition}
\label{prop:3.1}
Consider the numerical solution given by (\ref{AP-final}). Then,
\begin{itemize}
\item[$(i)$]  For all $\varepsilon \to 0$  and $\Delta t> 0$, the distribution function $f^{n+1}$ satisfies
$$
0 \leq f^{n+1}(x,v) \leq \max\left(\|f^n\|_{\infty}, \|\G[f^{n+1}] \|_{\infty}\right)
$$ 
\item[$(ii)$]  For all $\Delta t> 0$ and $f^0$, the distribution function $f^n$ converges to $\M[f^n]$, that is, 
$$
\lim_{\varepsilon\rightarrow 0} f^n = \M[f^{n}]
$$
and the  scheme gives a first order approximation in time 
of the compressible Euler system. 

\item[$(iii)$] Moreover, if we asssume that $\|f^n - \M^n\| = O(\var)$,  for $n\geq 2$ and
\begin{equation}
\label{hyp:0}
\left\|\frac{U^{n+1}-U^n}{\Delta t}\right\| \,\leq \, C,
\end{equation}
the scheme (\ref{AP-1})  asymptotically becomes a first order in time approximation of the compressible Navier-Stokes (\ref{eq:CNS}) given by
$$
\left\{
\begin{array}{l}
\displaystyle{\frac{\rho^{n+1} - \rho^n}{\Delta t} \,+\, \,{\rm div}_x\left( \rho^n \,u^n\right) \,=\, 0,}
\\
\,
\\
\displaystyle{\frac{\rho^{n+1} u^{n+1} - \rho^n u^n}{\Delta t}  \,+ \,{\rm div}_x \left(\rho^n \,u^n\otimes u^n + \rho\,T {\rm I} \right) \,=\, \var \,{\rm div}_x\left(\mu \sigma(u^{n-1})\right),}
\\
\,
\\
\displaystyle{\frac{E^{n+1} -  E^n}{\Delta t} \,+\,  {\rm div}_x \left(\left[E^n + \rho^n T^n\right]\,u^n\right)\,=\, \var\,{\rm div}_x\left(\,\mu \sigma(u^{n-1})\, u^{n-1} \,+\, \kappa \nabla_x T^{n-1}\right). }
\end{array}\right.
$$
where $\tau^n\,\mu\,=\, p^{n-1}/(1-\nu)$ and $\tau^n\,\kappa\,=\, (d_v+2)p^{n-1}/2$.
\end{itemize}
\end{proposition}
\begin{proof}
$(i)$ Let us first observe that $f^{n+1}$ is a linear combination of $f^n$ 
(defined along characteristics) and $\G[f^{n+1}]$, thus we get the first assertion.

To prove $(ii)$, for any initial distribution function, we considet $f^n$ for $n\geq 1$ and compute the asymptotic limit of $\Sigma^n$ when $\var$ goes to zero in (\ref{AP-final}), it yields
$$
\Sigma^n \,=\, \rho^n \,\left( T^n {\rm I} \,+\,  u^n\otimes u^n\right)  
$$ 
and using (\ref{sigm}), we also get
$$
\Theta^n  = T^n \, {\rm I},
$$
thus
$$
 \Tau^n = T^n \, {\rm I}\,.
$$
Hence in the asymptotic limit $\varepsilon\rightarrow 0$, the distribution function $\G[f^n]$ becomes isotropic, which means that $\G[f^n]$ converges
 to the classical Maxwellian $\M[f^n]$.

Therefore, the solution at zeroth order is obtained by taking $f^n=\M[f^n]$ in the conservation laws (\ref{AP-11}), namely,
$$
U^{n+1} =  \int_{\R^{d_v}}  \phi(v)\,(\M[f^n]-\Delta t v \cdot \nabla_x \M[f^n])  \, dv,
$$
and the scheme for the macroscopic quantities reduces to the well-known approximation to the Euler equation
\begin{equation}
\label{euler:approx}
\left\{
\begin{array}{l}
\displaystyle{\frac{\rho^{n+1} - \rho^n}{\Delta t} \,+  \,{\rm div}_x\left( \rho^n \,u^n\right)= 0,}
\\
\,
\\
\displaystyle{\frac{\rho^{n+1} u^{n+1} - \rho^n u^n}{\Delta t}  \,+\,{\rm div}_x \left(\rho^n \,u^n\otimes u^n + \rho^n\,T^n \,{\rm I}\right)= 0,  }
\\
\,
\\
\displaystyle{\frac{E^{n+1} -  E^n}{\Delta t} \,+\,  {\rm div}_x \left(\left[E^n + \rho^n T^n\right]\,u^n\right) = 0. }
\end{array}
\right.
\end{equation}  
Moreover, the temperature $T^{n+1}$ satisfies the following 
$$
\frac{d_v}{2}\,\frac{T^{n+1} - T^n}{\Delta t} + \frac{d_v}{2}\, u^n \cdot\nabla_x T^n  + \rho^n\,T^n\, {\rm div}_x u^n \,=\, O(\Delta t).  
$$
Now let us prove $(iii)$. The preservation of the asymptotic, that is the compressible Navier-Stokes equation,  is based on the Chapman-Enskog method. It simply consists of expanding the distribution function $f^n$ into 
$$
f^n = \M[f^n] \,+\, \varepsilon \,f_{1}^n, 
$$
which implies that 
\begin{equation}
\label{comp:cond}
\int_{\R^{d_v}} f_{1}^n(x,v)\,dv\,=\,0, \quad \int_{\R^{d_v}} f_{1}^n(x,v)\,v\,dv\,=\,0,\quad\int_{\R^{d_v}} f_{1}^n(x,v)\,|v|^2\,dv\,=\,0.
\end{equation} 
These conditions are known as the compatibility conditions. Moreover, we also expand the stress tensor  
\begin{equation}
\label{exp:Tau}
\Theta^n \,\,=\,\, T^n\, {\rm I} \,\,+\,\, \varepsilon \,\Theta_1^n
\end{equation}
and the heat flux  
\begin{equation}
\label{exp:Q}
\QQ^n(x) \,\,:=\,\, \int_{\R^{d_v}} \frac{|v-u^n|^2}{2} \,(v-u^n)\,f^n(x,v)dv  \,\,=\,\,  0 \,\,+\,\, \varepsilon\; \QQ_1^n(x),
\end{equation}
where 
$$
\Theta_1^n \,=\, \frac{1}{\rho^n}\int_{\R^{d_v}} f_1^n(x,v)(v-u^n)\otimes (v-u^n)\,dv, \quad \QQ_1^n \,=\, \int_{\R^{d_v}} f_1^n(x,v)\,\frac{|v-u^n|^2}{2} \, (v-u^n)\,dv
$$
and ${\rm tr}\left(\Theta_1^n\right)=0$. Inserting this latter expansion into the discrete conservation laws (\ref{AP-11}), it gives
$$
\left\{
\begin{array}{l}
\displaystyle{\frac{\rho^{n+1} - \rho^n}{\Delta t} \,+\, \,{\rm div}_x\left( \rho^n \,u^n\right) \,=\, 0,}
\\
\,
\\
\displaystyle{\frac{\rho^{n+1} u^{n+1} - \rho^n u^n}{\Delta t}  \,+ \,{\rm div}_x \left(\rho^n \,u^n\otimes u^n + \rho^n\,T^n{\rm I}\right) \,=\, -\var \,{\rm div}_x\left(\rho^n\,\Theta_1^n\right),}
\\
\,
\\
\displaystyle{\frac{E^{n+1} -  E^n}{\Delta t} \,+\,  {\rm div}_x \left(\left[E^n + \rho^n T^n\right]\,u^n\right)\,=\, -\var\,{\rm div}_x\left(\,\QQ_1^n\,+\,{\rho^n\,\Theta_1^n}\,u\,\right). }
\end{array}\right.
$$
For the application of the Chapman-Enskog method, the anisotropic Gaussian $\G(f)$ must be expanded with respect to $\var$ as well,
$$
\G[f^n] \,=\, \M[f^n] \,\,+\,\, \var\, g_1^n
$$
The next step is to insert these expansions into the scheme  (\ref{AP-1}) for the ES-BGK eqution  and  use the compatibility conditions, and it yields for $n\geq 1$
\begin{eqnarray}
\label{4.19}
 \frac{\M[f^{n}]-\M[f^{n-1}]}{\Delta t} \,+ \,v\cdot\nabla_x \M[f^{n-1}] &=& \tau^n\,(g_1^{n} - f_1^{n})  
\\
&-& \var\left[\frac{f_1^{n}-f_1^{n-1}}{\Delta t}  + \,v\cdot\nabla_x f^{n-1}_1\right].
\nonumber
\end{eqnarray}
One the one hand, the term $g_1^n$ is computed from $\G[f^n]$ as follows. First, the distribution function $\G[f^n]$ contains the inverse matrix $\Tau^n$ for which insertion of (\ref{exp:Tau}) yields
$$
\Tau^n  \,\,=\,\, T^n \,{\rm I} \,\,+\,\, \nu\,\varepsilon\,\Theta_1^n.
$$
Therefore, a Taylor expansion within terms of second order and using the fact that ${\rm tr}(\Theta_1^n)=0$,  gives on the one hand
$$
{\rm det}(\Tau^n) \,\,=\,\, (T^n)^{d_v} \,\,+\,\, O(\varepsilon^2) 
$$
and on the other hand
$$
[\Tau^n]^{-1} \,\,=\,\, \frac{1}{T^n} \,\left[{\rm I} \,-\,\frac{\nu \,\varepsilon}{T^n} \Theta_1^n \right]  \,+\, O(\varepsilon^2) 
$$
and therefore
$$
g_1^n \,\,=\,\, \frac{\G[f^n]-\M[f^n]}{\var} \, \,= \,\,  \nu\,\frac{\M[f^n]}{2\,(T^n)^2} \,(v-u^n)^t\,\Theta_1^n \,(v-u^n).
$$
On the other hand, to obtain the first order expression for the distribution function $f_1^n$, we need to consider the terms of order zero in (\ref{4.19}), which we write for convenience as
$$
f_1^n \,=\, g_1^n -  \frac{\M[f^{n}]-\M[f^{n-1}]}{\tau^n\,\Delta t} \,+ \,\frac{v}{\tau^n}\cdot\nabla_x \M[f^{n-1}] \,\,+\,\, O(\var).
$$
The differential ${\rm d}\M[f]$ of the Maxwellian is given by
$$
{\rm d}\M[f] \,\,=\,\, \M[f]\,{\rm d}\log(\M[f]) \,\,=\,\, \M[f] \left[\frac{{\rm d}\rho}{\rho} + \left(\frac{|v-u|^2}{2 T} - \frac{d_v}{2}\right)\frac{{\rm d}T}{T} - \frac{u-v}{T}{\rm d}u \right]
$$
and 
with the assumption (\ref{hyp:0}), we obtain
\begin{eqnarray*}
\frac{\M[f^{n}]-\M[f^{n-1}]}{\Delta t} \,\,=\,\,\,\, \M[f^{n-1}] &&\left[\frac{1}{\rho^{n-1}}\,\frac{\rho^{n}-\rho^{n-1}}{\Delta t}  \right. + \frac{1}{T^{n-1}}\,\left(\frac{|v-u^{n-1}|^2}{2 T^{n-1}}  - \frac{d_v}{2}\right)\frac{T^{n}-T^{n-1}}{\Delta t}
\\
&&- \left.\frac{u^{n-1}-v}{T^{n-1}} \cdot\frac{u^{n}-u^{n-1}}{\Delta t} \right] \,+\, O(\Delta t) 
\end{eqnarray*}
and
$$
v\,\cdot\,\nabla_x \M[f^{n-1}] =  \M[f^{n-1}]\left[\frac{v\cdot\nabla_x\rho^{n-1}}{\rho^{n-1}} + \left(\frac{|v-u^{n-1}|^2}{2 T^{n-1}} - \frac{d_v}{2}\right)\frac{v\cdot\nabla_xT^{n-1}}{T^{n-1}} - \frac{u^{n-1}-v}{T^{n-1}} v\cdot\nabla_x u^{n-1} \right].
$$
Gathering the two latter equalities, it yields up to the order $\Delta t$
\begin{eqnarray*}
\lefteqn{\frac{\M[f^{n}]-\M[f^{n-1}]}{\Delta t} \,+\, v\cdot\nabla_x \M[f^{n-1}] \,=\, } 
\\
&&  \frac{\M[f^{n-1}]}{T^{n-1}} \left[(u^{n-1}-v)\,\nabla_x u^{n-1} \,(u^{n-1}-v) + (u^{n-1}-v)\left(\frac{|v-u^{n-1}|^2}{2 T^{n-1}} - \frac{d_v+2}{2}\right)\nabla_x T^{n-1} \right]
\\
&+& \frac{\M[f^{n-1}]}{\rho^{n-1}} \left[ \frac{\rho^{n}-\rho^{n-1}}{\Delta t} + {\rm div}_x(\rho^{n-1}\,u^{n-1}) \right] 
\\
&+& \frac{\M[f^{n-1}]\,(v-u^{n-1})}{T^{n-1}} \left[ \frac{u^{n}-u^{n-1}}{\Delta t} +  u^{n-1}\cdot\nabla_x u^{n-1} +  \frac{1}{\rho^{n-1}} \nabla_x (\rho^{n-1}\,T^{n-1}) \right]
\\
&+& \M[f^{n-1}]\, \left(\frac{|v-u^{n-1}|^2}{2 T^{n-1}} - \frac{d_v}{2}\right) \,\left[ \frac{1}{T^{n-1}}\left(\frac{T^{n}-T^{n-1}}{\Delta t} +  u^{n-1}\cdot\nabla_x T^{n-1}\right)  \,+\,  \frac{2}{d_v} {\rm div}_x u^{n-1} \right].   
\end{eqnarray*}
The last three lines may vanish due to the conservation laws (\ref{euler:approx}) up to the order of $\var$ and $\Delta t$. Thus, the result for the first order contribution to the distribution function is
\begin{eqnarray*}
f_1^n &=& -\frac{\M[f^{n-1}]}{\tau^n\,T^{n-1}} \left[(u^{n-1}-v)\,\nabla_x u^{n-1} \,(u^{n-1}-v) + (u^{n-1}-v)\left(\frac{|v-u^{n-1}|^2}{2 T^{n-1}} - \frac{d_v+2}{2}\right)\nabla_x T^{n-1} \right]
\\
&+&  \nu\,\frac{\M[f^n]}{2\,|T^n|^2} \,(v-u^n)\,\Theta_1^n \,(v-u^n)  \,+\, O(\Delta t) \,+\, O(\var).
\end{eqnarray*}
It is straightforward to show that $f_1$ satisifes the compatibility relations (\ref{comp:cond}). Hence, we get for the stress tensor $\Theta_1^n$
$$
\rho^n\,\Theta_1^n \,=\, \int_{\R^{d_v}}f_1^n (v-u^n)\otimes (v-u^n)dv \,\,=\,\, -  \frac{\rho^{n-1}\, T^{n-1}}{\tau^n}\,\sigma(u^{n-1})  \,\,+\,\,\nu \,\rho^n \,\Theta_1^n  \,\,+\,\, O(\Delta t) \,+\, O(\var),
$$
that is,
$$
\rho^n\,\Theta_1^n  \,\,=\,\, -  \frac{p^{n-1}}{(1-\nu)\tau^n}\,\sigma(u^{n-1})  \,+\, O(\Delta t) \,+\, O(\var).
$$
Then the heat flux $\QQ_1^n$ is given by
$$
\QQ_1^n \,=\, \int_{\R^{d_v}} \frac{|v-u^n|^2}{2} \,(v-u^n)\,f^n_1\,dv \,=\,  - \frac{d_v+2}{2} \,\frac{p^{n-1}}{\tau^n}\,\nabla_x T^{n-1}\,+\, O(\Delta t) \,+\, O(\var).
$$ 

The constitutive relations are laws of Navier-Stokes and Fourier where
$$
\mu \,=\, \frac{1}{1-\nu}\,\frac{p^{n-1}}{\tau^n}\qquad{\rm and}\qquad\kappa\,=\, \frac{d_v+2}{2} \,\frac{p^{n-1}}{\tau^n} 
$$
are viscosity and thermal conductivity, respectively. Thus, a first order approximation with respect to $\Delta t$ and $\var$ is given by
$$
\rho^n\,\Theta_1^n  \,\,=\,\, -  \mu\, \sigma(u^{n-1})  \qquad{\rm and}\qquad \QQ_1^n \,=\,-\kappa \nabla_x T^{n-1}
$$
The Prandtl number is related to the coefficient $\nu$ of the ES-BGK model by
$$
{\rm Pr}  \,=\, \frac{d_v+2}{2}\frac{\mu}{\kappa} = \frac{1}{1-\nu}. 
$$  
\end{proof}

\section{Numerical simulations}

In this section, we give three numerical examples for the ES-BGK equation in different asymptotic regimes
in order to check the performance (in stability and
accuracy) of our methods. 
We have implemented the first order  scheme  (\ref{AP-final}) for the approximation of the ES-BGK equation. A classical second 
order finite volume scheme with slope limiters is applied for the transport 
operator. We present two numerical tests for a $1d_x\times 2d_v$ model and  finally a non stationary $2d_x\times 2d_v$ model.

We will compare the numerical solution to the ES-BGK equations with the one obtained for the full Boltzmann equation with Maxwellian molecules using a spectral approximation \cite{FMP, FJ} and the one obtained for the compressible Navier-Stokes system using a second order finite volume scheme.

To this aim, we need to choose the right value for $\tau$ such that the viscosity and heat conductivity computed from the asymptotic limit of the ES-BGK model is the same as the ones corresponding to the full Boltzmann equation. According to \cite{chapman}, the viscosity computed for the full Boltzmann equation is given by
$$
\mu_B(T) = \frac{\sqrt{2}}{3\pi}\,\frac{T}{A_2(5)},
$$ 
where $A_2(5)\simeq 0.436$. On the other hand, the viscosity computed from the ES-BGK model is
$$
\mu = \frac{1}{1-\nu} \frac{p}{\tau},
$$ 
where $\nu=-1$. Thus, we choose $\tau$ such that both viscosity are equal, which leads to
$$
\tau = \frac{1}{2} \frac{p}{\mu_B(T)} = \frac{3\pi}{2\sqrt{2}}\,A_2(5)\rho \simeq 0.925\,\frac{\pi}{2}\,\rho.
$$

\subsection{Approximation of smooth solutions}
For this numerical test, we consider the ES-BGK equation in dimension $1d_x\times 2d_v$ on the torus
$$
\left\{
\begin{array}{l}
\displaystyle{\frac{\partial f}{\partial t} + v \cdot \nabla_x f = \frac{1}{\var}Q(f,f), \quad x\in [-1,1], v \in \R^2}
\\
\,
\\
f(t=0)=f_0,
\end{array}
\right.
$$
with   periodic boundary conditions in $x$. The operator $\Q(f)$ is given by $\Q(f) \,=\, \tau\,[\G[f]-f]$ where the parameter  $\tau$ is chosen in order that the viscosity $\mu$ matches perfectly with one obtained to the full Boltzmann operator for Maxwellian molecules, that is $\tau = 0.9\,\pi\,\rho/2$. 

Define $\rho_g$ and $T_g$ with respect to the initial data $f_0$ by    $$
\rho_g = \frac{1}{2}\int_{-1}^1 \int_{\R^2} f_0(x,v) dv dx,\quad{\rm and}\quad   T_g = \frac{1}{2\rho_g}\int_{-1}^1 \int_{\R^2} f_0(x,v)\,|v|^2\, dv dx
$$
and assume for simplicity  that
$$
\frac{1}{2}\int_{-1}^1 \int_{\R^2} f_0(x,v)\,v\, dv dx= 0.  
$$
Whenever $f(t,x,v)$ is a smooth solution to the Boltzmann or ES-BGK equation with periodic boundary conditions, one has the global conservation laws for mass, momentum and energy. These conservation laws are then enough to uniquely determine the stationary state of the model : the normalized global Maxwellian distribution
\begin{equation}
\label{maxwglob}
\M_g(v) = \frac{\rho_g}{2\pi\,T_g}\,\exp\left(-\frac{|v|^2}{2T_g}\right),\quad v\in\R^2.
\end{equation}

Our goal here is to investigate numerically the long-time behavior of the solution $f$ and to compare the solution with the asymptotic behavior of the solution to the compressible Navier-Stokes equations (CNS).  If $f$ is any reasonable solution of the ES-BGK equation, satisfying certain {\em a priori} bounds of compactness (in particular, ensuring that no kinetic energy is allowed to leak at large velocities), then it is expected that $f$ does indeed converge to the global Maxwellian distribution $\M_g$ as $t$ goes to $+\infty$. 

Recently, Desvillettes and Villani \cite{DV:EB:03}, Guo and Strain \cite{Guo:rate} were interested in the study of rates of convergence for the full Boltzmann equation. Roughly speaking in \cite{DV:EB:03}, the authors proved that if the solution to the Boltzmann equation is smooth enough  then (with constructive bounds)
$$
\| f(t) - \M_g \| = O(t^{-\infty}),
$$
which means that the solution converges almost exponentially fast to the global equilibrium (namely with polynomial rate $O(t^{-r})$ with $r$ as large as wanted). Moreover in \cite{DV:EB:03}, Desvillettes and Villani conjectured that time oscillations should occur on the evolution of the relative local entropy
$$
\HH_l(t) = \int f \, \log\left(\frac{f}{\M_l}\right) \, dx \, dv,
$$
where  $\M_l$ is the local Maxwellian distribution in the sense that the constants
$\rho$, $u$ and $T$ appearing there depend on  time $t$ and position $x$
\begin{equation}
\label{maxwloc}
\M_l(t,x,v) = \frac{\rho(t,x)}{2\pi T(t,x)}\,\exp\left(-\frac{|v-u(t,x)|^2}{2T(t,x)}\right)
\end{equation}
 In fact their proof does not rule out the possibility that the entropy production undergoes important oscillations in time, and actually most of the technical work is caused by this possibility.

To estimate the speed of convergence to the global equilibrium and the possibility that oscillations also occur on the difference between global and local equilibria (since the global relatve entropy is nonincreasing), we prefer to investigate the behavior of the following quantity
$$
\mathcal{E}(t) := \|\rho(t) - \rho_g \|_{L^1},
$$ 
which makes sense both for the solution to the Boltzmann or ES-BGK equation but also for the compressible Navier-Stokes equation.

Here, we performed simulations to the full Boltzmann equation  with a fast spectral method \cite{FMP}, to the ES-BGK equation with the scheme (\ref{AP-1}) and to the compressible Navier-Stokes system with a WENO solver in a simplified geometry (one dimension in space, two dimension in velocity, periodic boundary conditions) to  check numerically if such oscillations occur. Clearly this test is challenging for a numerical method due to the high accuracy required to capture such oscillating behavior at the kinetic regime.

Then, we consider an initial datum as a perturbation of the global equilibrium $\M_g$
\begin{equation*}
f_0(x,v) = \frac{1+ A_0 \, \sin(\pi\,x)}{2\pi T_0}\, \left[\exp\left(-\frac{|v-u_0|^2}{2T_0}\right)+\exp\left(-\frac{|v+u_0|^2}{2T_0}\right)\right], x \in [-1,1], v \in \R^2,
\end{equation*}
with $A_0=0.5$, $T_0=0.125$ and $u_0=(1/2,1/2)$. We have chosen $\nu=-1$ such that the Prandtl number of the ES-BGK model corresponds to the Prandtl number of the $2d_v$ Boltzmann operator, that is ${\rm Pr}=0.5$.

\begin{figure}[htbp]
\begin{tabular}{cc}
\includegraphics[width=7.cm]{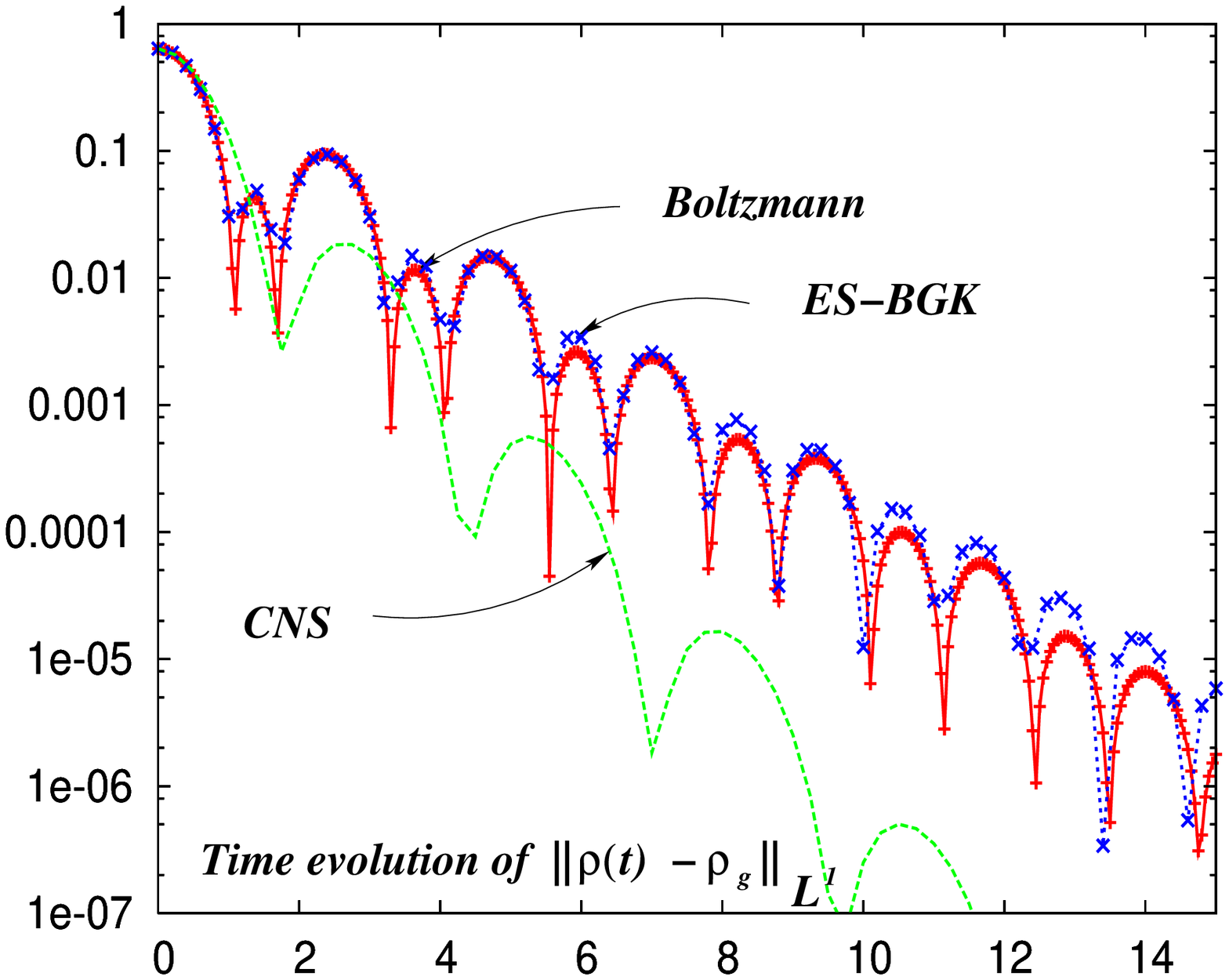}  &
\includegraphics[width=7.cm]{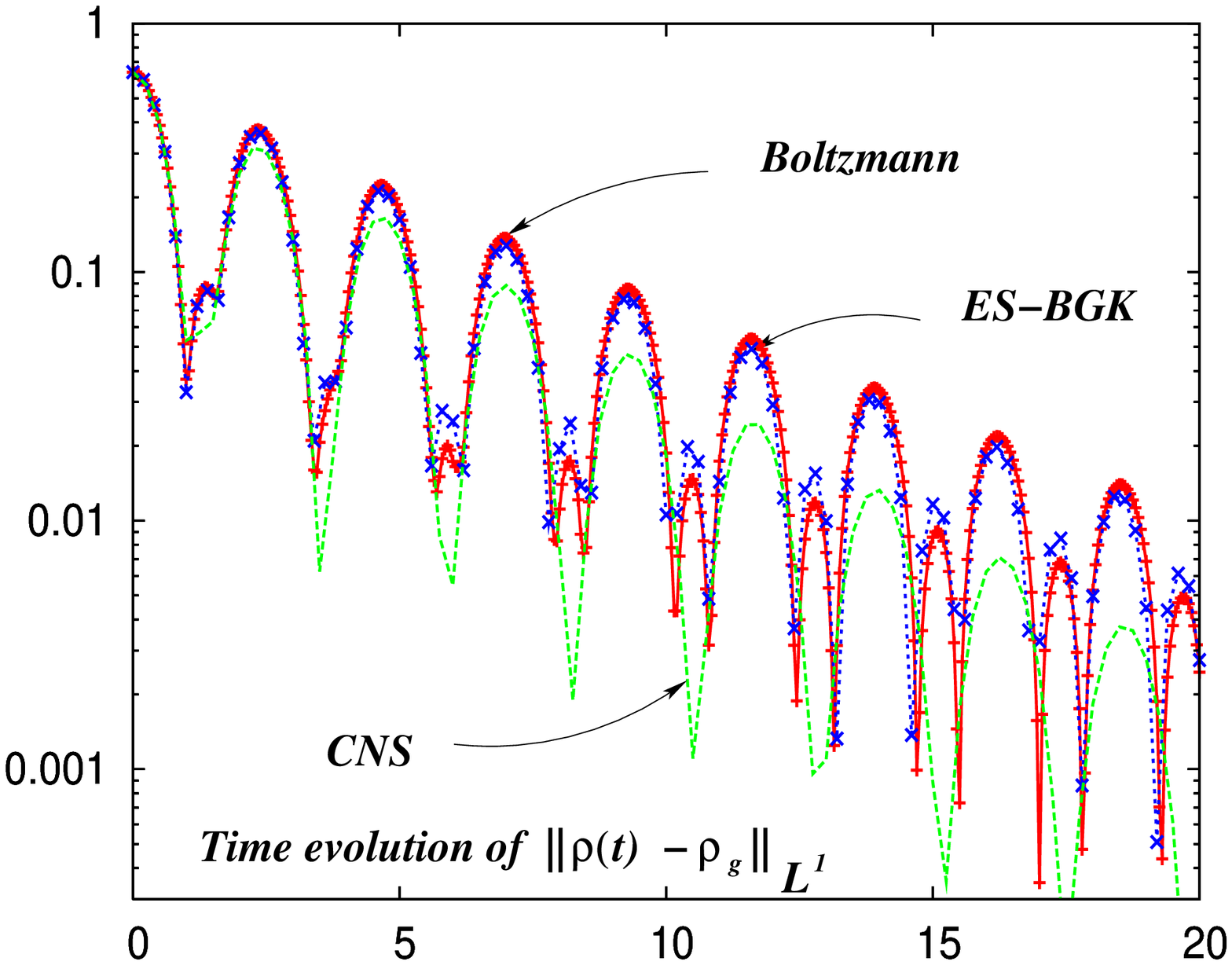} 
\\
(a) & (b)
\\
\includegraphics[width=7.cm]{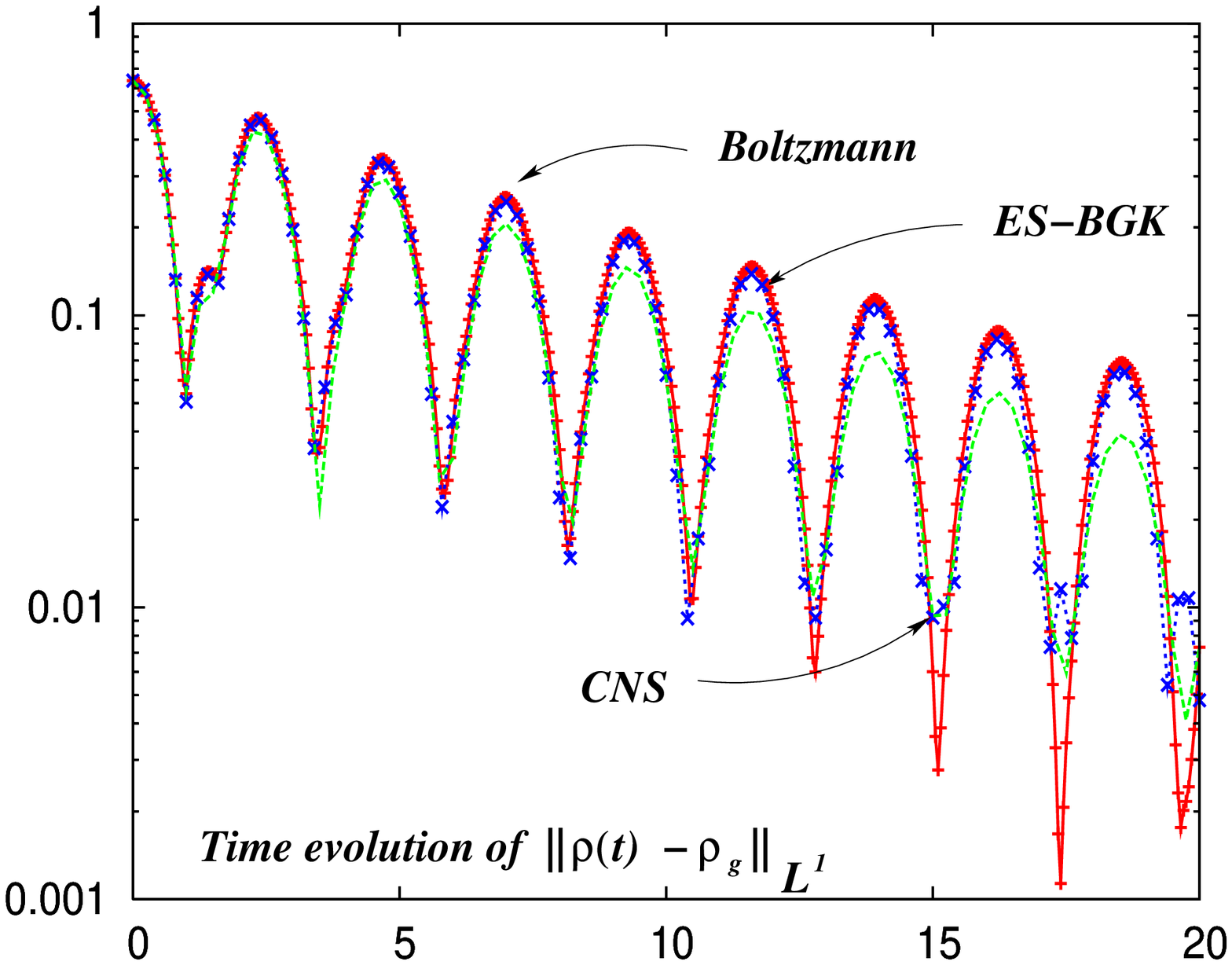}& 
\includegraphics[width=7.cm]{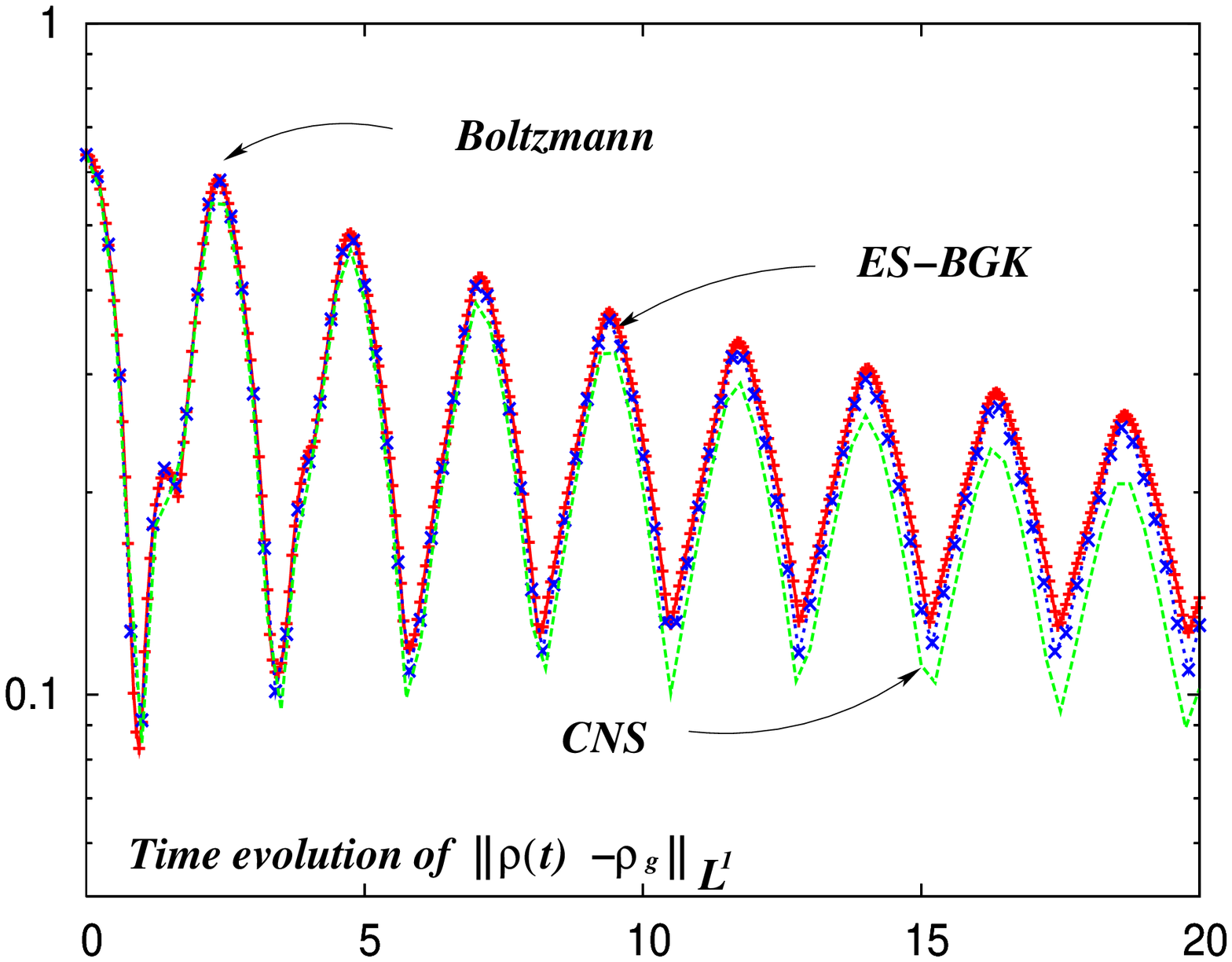} 
\\
(c) & (d) 
\end{tabular}
\caption{Influence of the Knudsen number $\var$: distance between the local  density $\rho(t)$ and the global density at equilibrium $\rho_g$  using $100 \times 32\times 32$ for (a) $\var$=$0.5$; (b) $\var=0.1$; (c) $\var=0.05$ and (d) $\var=0.01$.} 
\label{fig3-1}
\end{figure}

In Figure \ref{fig3-1}, one can indeed  observe oscillations on the quantitiy $\mathcal{E}(t)$ for the full Boltzmann equation, the ES-BGK model and also for the compressible Navier-Stokes system obtained from the asymptotic of the ES-BGK equation (\ref{eq:esbgk}). The strength of the oscillations does not depend  on the Knudsen number $\var$. The superimposed curves yield the time evolution of the $\mathcal{E}(t)$ for $t\in [0,20]$; the first plot corresponds to $\var=0.5$, the second one to $\var=0.1$, the third one $\var=0.05$ and the last one $\var=0.01$. 

On the one hand, for $\var=0.5$, which corresponds to a rarefied regime, the behavior of $\mathcal{E}(t)$ strongly differs between the kinetic and hydrodynamic models. The results for ES-BGK and the full Boltzmann equations agree very well in this rarefied regime, which illustrates perfectly the efficiency and accuracy  of the ES-BGK model.

On the other hand, for smaller values of  $\var$, the different numerical approximations give roughly the same results and the ES-BGK model and  the compressible Navier-Stokes system become very close. Finally for $\var\simeq 0.01$, the different kinetic  (Boltzmann and ES-BGK) and hydrodynamic (Euler or Navier-Stokes) models agree very well. 

Further note that the equilibration is much more rapid when the Knudsen number $\var$ is large, and that the convergence seems to be exponential.

\subsection{The Riemann problem}
This test deals with the numerical solution to the $1d_x\times 2d_v$ ES-BGK equation. The operator $\Q(f)$ is given by $\Q(f) \,=\, \tau\,[\G[f]-f]$ where $\tau = 0.9\,\pi\rho/2$.  We present several numerical simulations for one dimensional Riemann problem, with different Knudsen numbers, from rarefied regime  to the fluid regime.

Here, the initial data corresponding to the ES-BGK equations are given by 
the isotropic Maxwellian distributions computed from the following macroscopic quantities
\begin{eqnarray*}
\left\{
\begin{array}{ll}
(\rho_l, u_l, T_l) = (1,M\,\sqrt{2},1)\,, &\textrm{ if } x \leq 0\,,
\\ 
\\
(\rho_r, u_r, T_r) = (1.,0,1.05)\,, &\textrm{ if } x > 0\,
\end{array} \right.
\end{eqnarray*}
with the Mach number $M=2.5$. We perform several computations for $\varepsilon\,=\, 5\times 10^{-1}$, $10^{-1}$,..., $10^{-3}$.

We present a comparison between the numerical solution to the Boltzmann equation obtained using a spectral scheme \cite{FMP}, the approximation to the compressible Navier-Stokes system obtained using a fifth order WENO  and our first order implicit method (\ref{AP-1}) for the ES-BGK model with $\nu=1/2$, for which the Prandtl number of the ES-BGK model is the same as the one for the $2d_v$ Boltzmann operator. Let us note that the viscosity and conductivity used for the numerical simumation of the compressible Navier-Stokes system are the ones obtained from the Chapman-Enskog expansion.

\begin{figure}[htbp]
\begin{tabular}{cc}
\includegraphics[width=7.5cm]{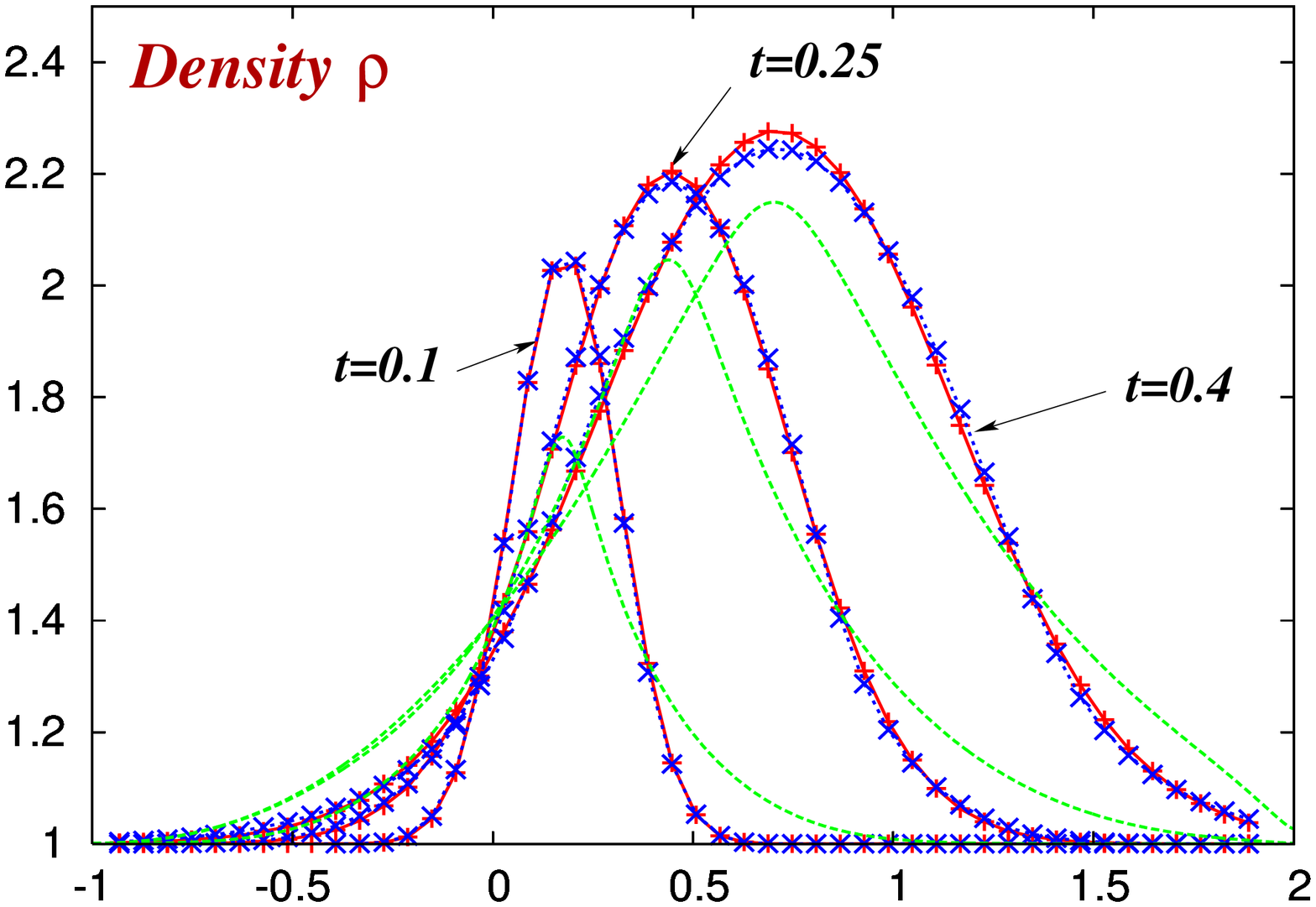}    
&
\includegraphics[width=7.5cm]{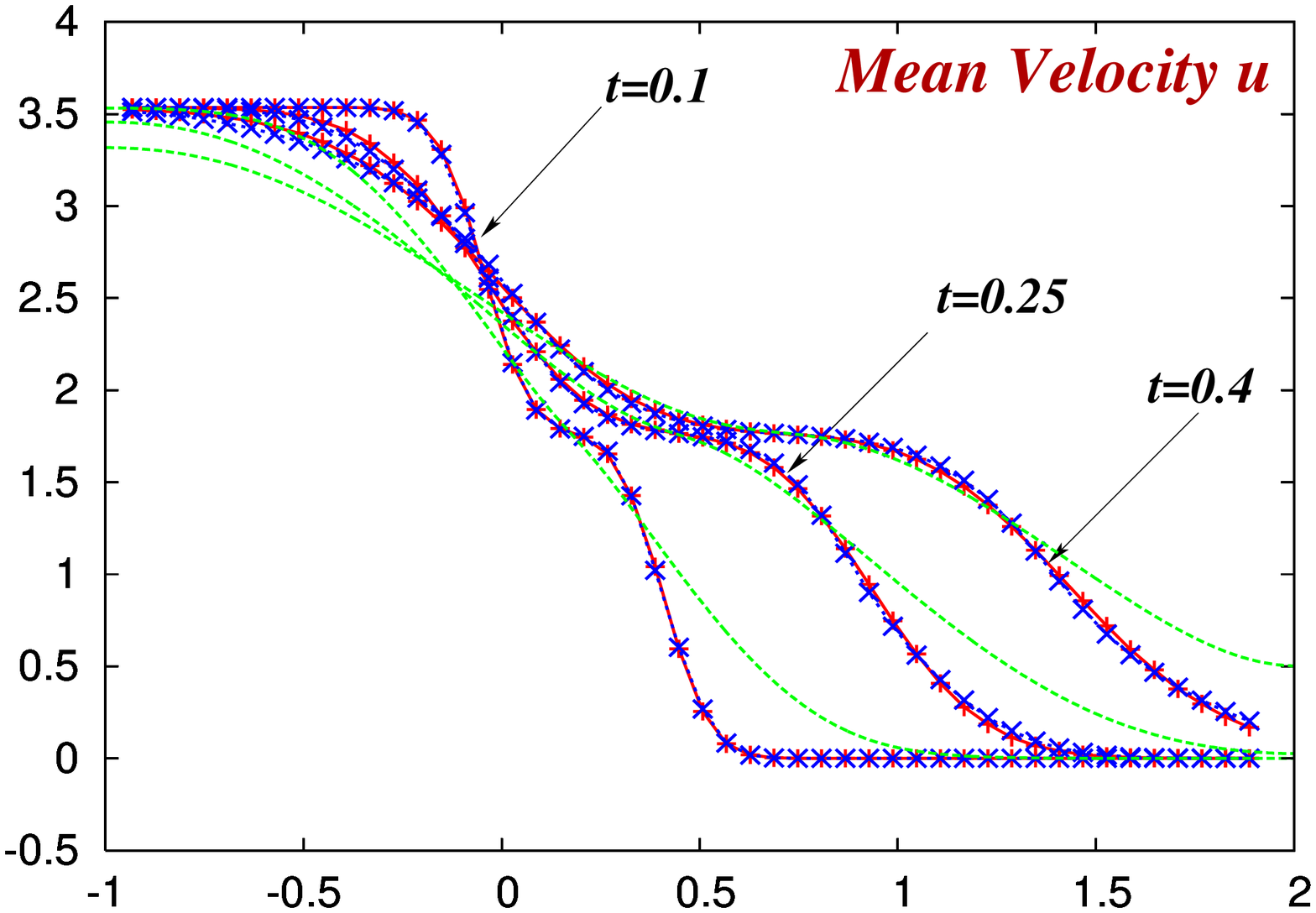}    
\\
(1)&(2)
\\
\includegraphics[width=7.5cm]{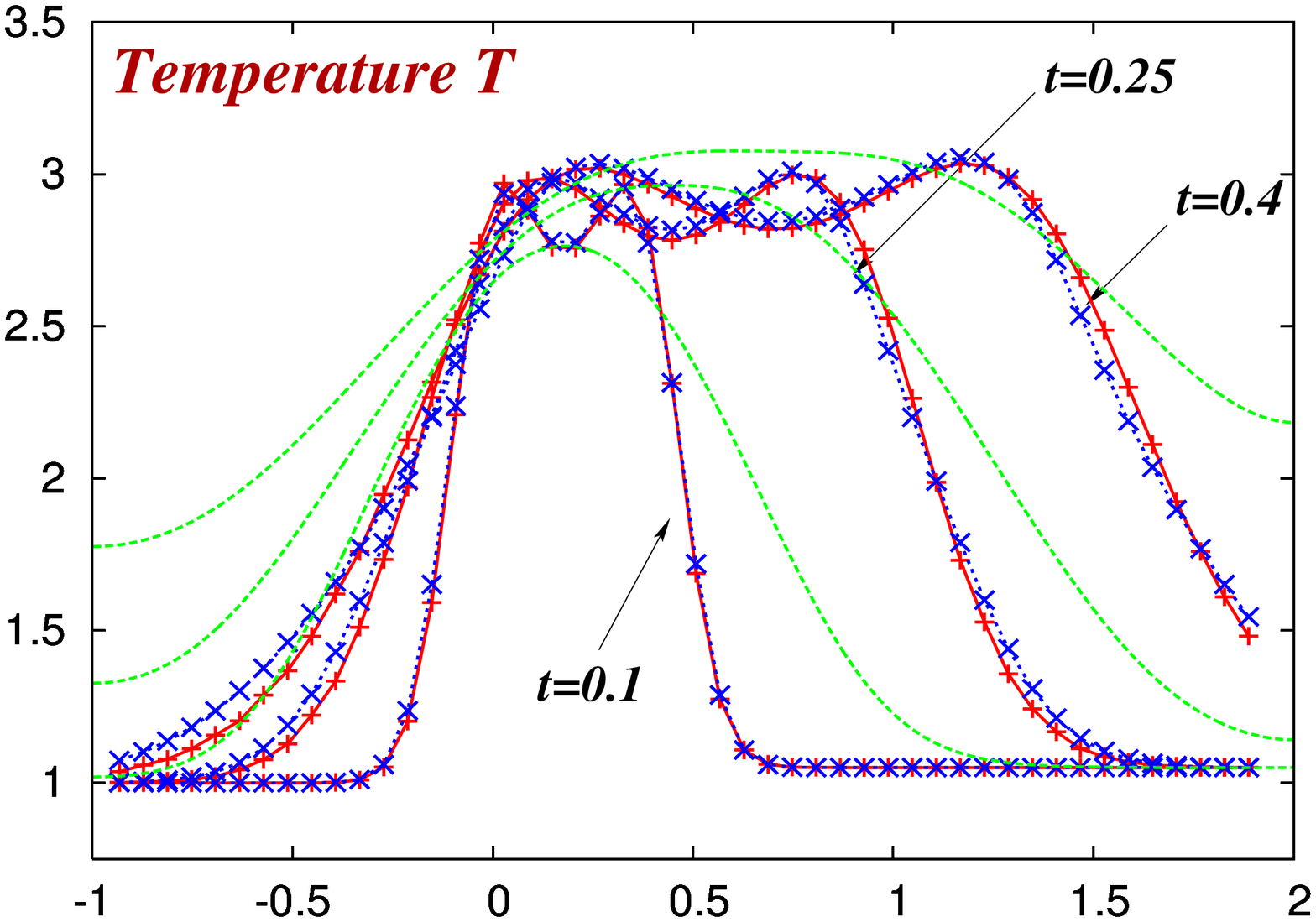}   
&
\includegraphics[width=7.5cm]{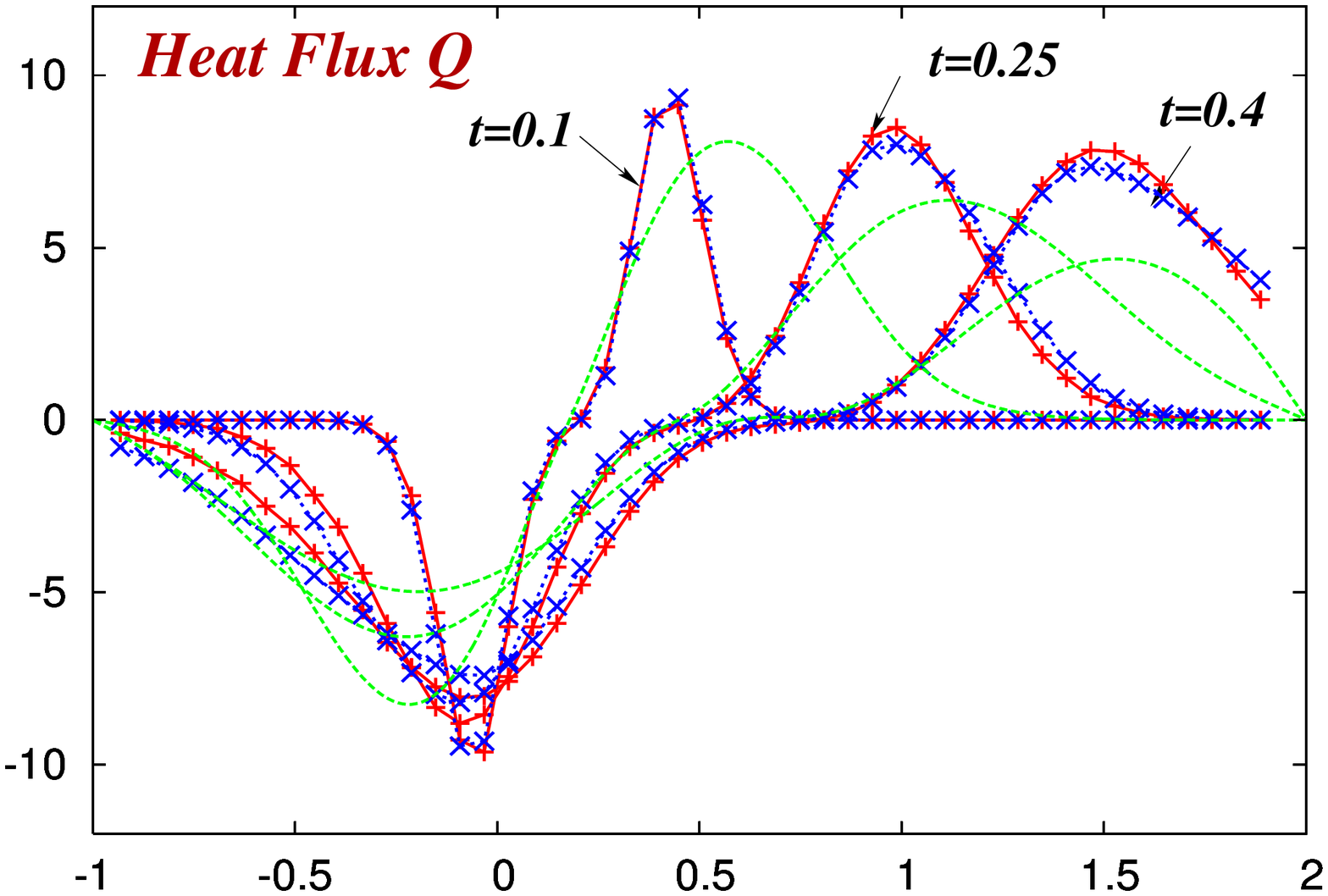}   
\\
(3)&(4)
\end{tabular}
\caption{Riemann problem ($\varepsilon=5\times 10^{-1}$):  crosses ({\tt +}) represent the numerical solution to the Boltzmann equation  obtained with our method (\ref{AP-1}), stars ({\tt x}) represent the solution  to the ES-BGK model  and lines is the solution corresponding to the compressible Navier-Stokes system. Evolution of (1) the density $\rho$, (2) mean velocity $u$, (3) temperature $T$ and  (4) heat flux $\QQ$ at time $t=0.1$, $0.25$ and $0.4$.}
\label{fig:03-2}
\end{figure}

In Figuer~\ref{fig:03-2}, we take $\var=5\times 10^{-1}$ and choose a time step $\Delta t=0.001$ satisfying the CFL condition for the transport part (with $n_x=200$). For such a value of $\varepsilon$, the problem is not stiff and this test is only performed to compare the accuracy of our scheme (\ref{AP-1}) with the different models. We present several snapshots of the density, mean velocity, temperature and heat flux 
$$
\QQ_1(t,x) \,:=\, \frac{1}{\varepsilon}\,\int_{\R^{d_v}}\frac{|v-u|^2}{2}\, (v-u) \,f(t,x,v)dv 
$$
at different time $t=0.1$, $0.25$ and $0.4$. We observe that, for a short time $t=0.1$, the numerical approximation of macroscopic quantities  and  heat flux given by (\ref{AP-1}) for the ES-BGK model are relatively close to the numerical solution to the Boltzmann equation. The front speed and the shape of the temperature with two bumps are very well approximated by the ES-BGK model, and the heat flux given by the Boltzmann equation and the ES-BGK model are different from the ones given by the compressible Navier-Stokes system. Clearly, the compressible Navier-Stokes system is not adequate in this rarefied regime whereas the ES-BGK model gives very satisfying results.

\begin{figure}[htbp]
\begin{tabular}{cc}
\includegraphics[width=7.5cm]{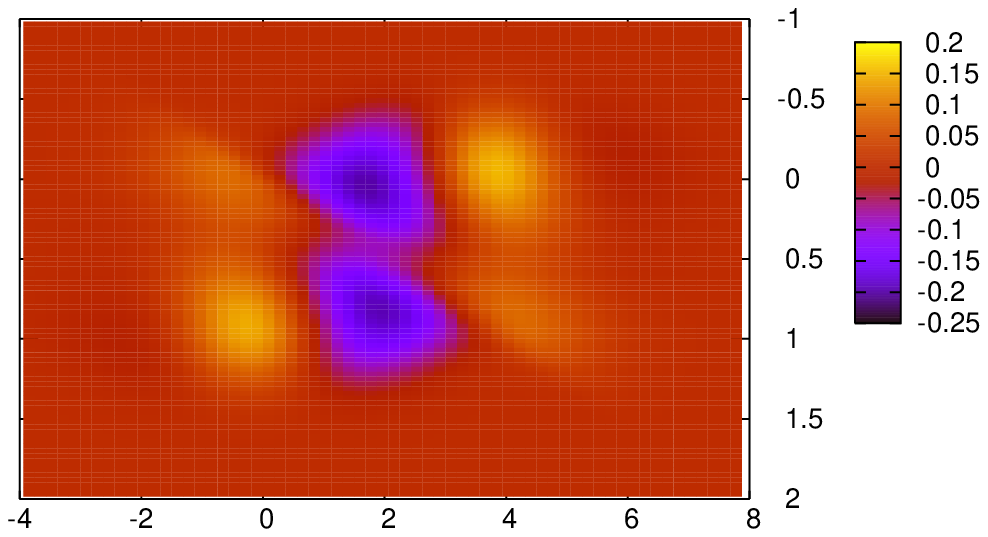}    
&
\includegraphics[width=7.5cm]{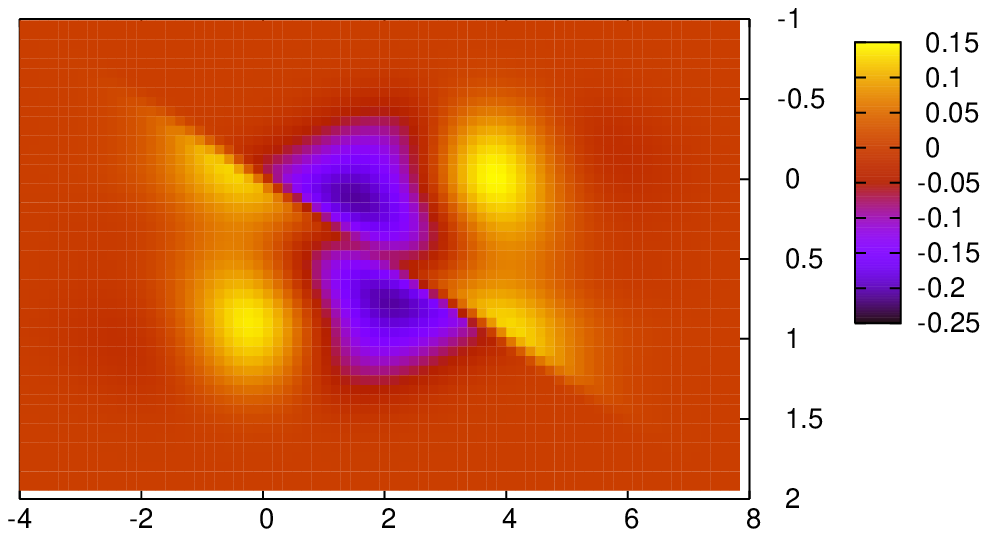}    
\\
(1)&(2)
\end{tabular}
\caption{Riemann problem ($\varepsilon=5\times 10^{-1}$):  $x-v_x$ projection of $f-\M[f]$ at time $0.25$ for the $(1)$ Boltzmann equation and $(2)$ ES-BGK model.}
\label{fig:f}
\end{figure}

For a larger time $t\geq 0.25$, the distribution in velocity $f$ to the ES-BGK model and to the Boltzmann equation agree very well and the macroscopic quantities are already well approximated by the solution to the compressible Navier-Stokes system. In Figure~\ref{fig:f}, we represent the $x-v_x$ projection of   $f-\M[f]$ at time $t=0.25$ obtained for the Boltzmann equation and the ES-BGK model. The distribution function becomes particularly far from the equilibrium at the front of the shock in velocity and then propagates in the computational domain. As can be observed on the macroscopic quantities, the solution to the Boltmzann equation is very close to the one obtained from the ES-BGK model.
    
Finally, at the kinetic regime our method (\ref{AP-1}) gives the same accuracy as a standard first order fully explicit scheme for the ES-BGK model or full Boltzmann equation without any additional computational effort. Of course, the computational effort needed for the ES-BGK models is much  smaller than the one for an accurate discretization of the full Boltzmann equation.

Now, we investigate the cases of small values of $\varepsilon$ for which an 
explicit scheme requires the time step to be of order $O(\var)$. In order to 
evaluate the accuracy of our method (\ref{AP-1}) in the Navier-Stokes regime (for small $\varepsilon\ll 1$ but not negligible), we  compared  the numerical solution  for $\varepsilon=10^{-1}$ with one obtained by the approximation of the compressible Navier-Stokes system derived from the ES-BGK model since the viscosity and heat conductivity are in that case explicitly known \cite{BLM}.

\begin{figure}[htbp]
\begin{tabular}{cc}
\includegraphics[width=7.5cm]{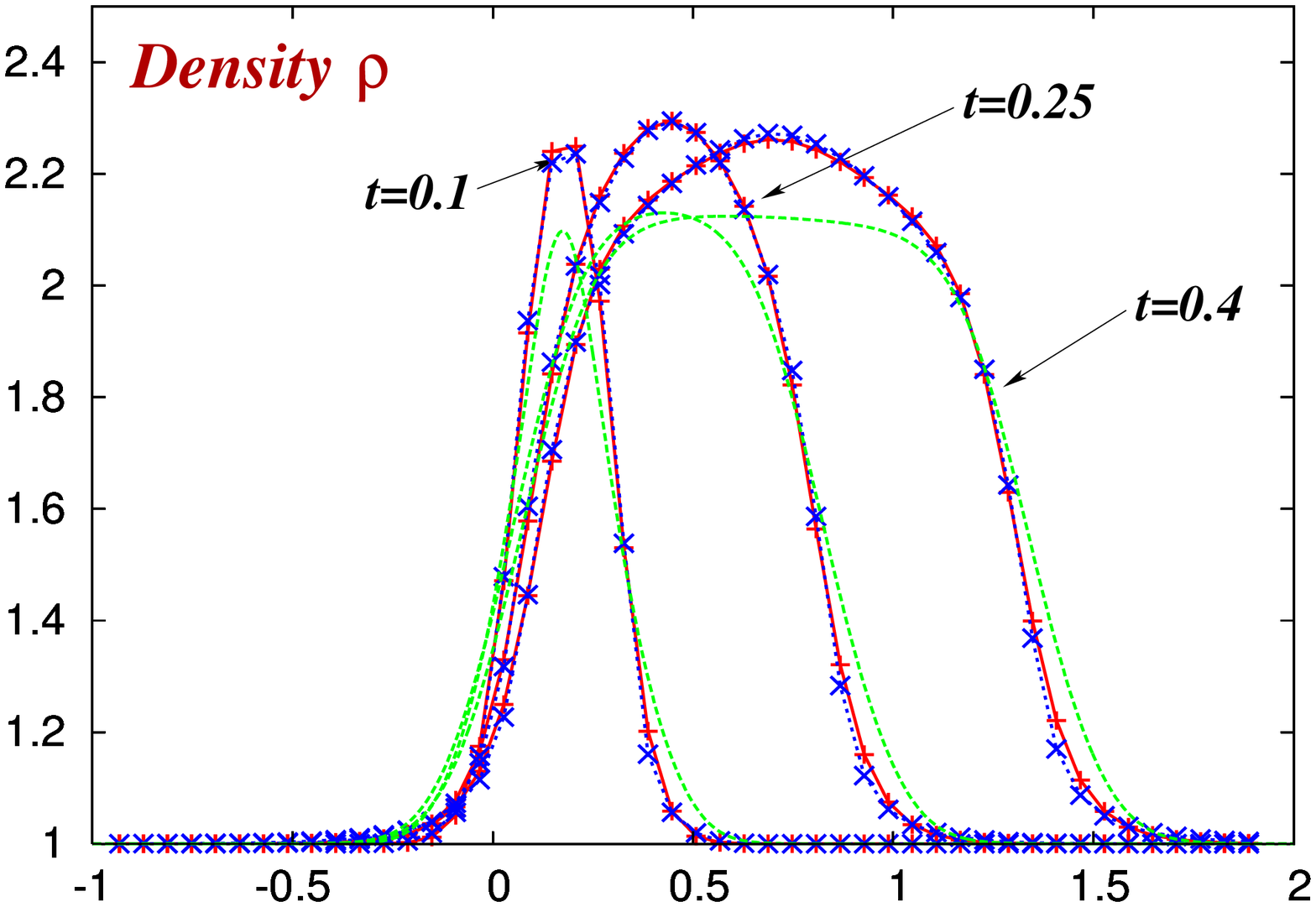}    
&
\includegraphics[width=7.5cm]{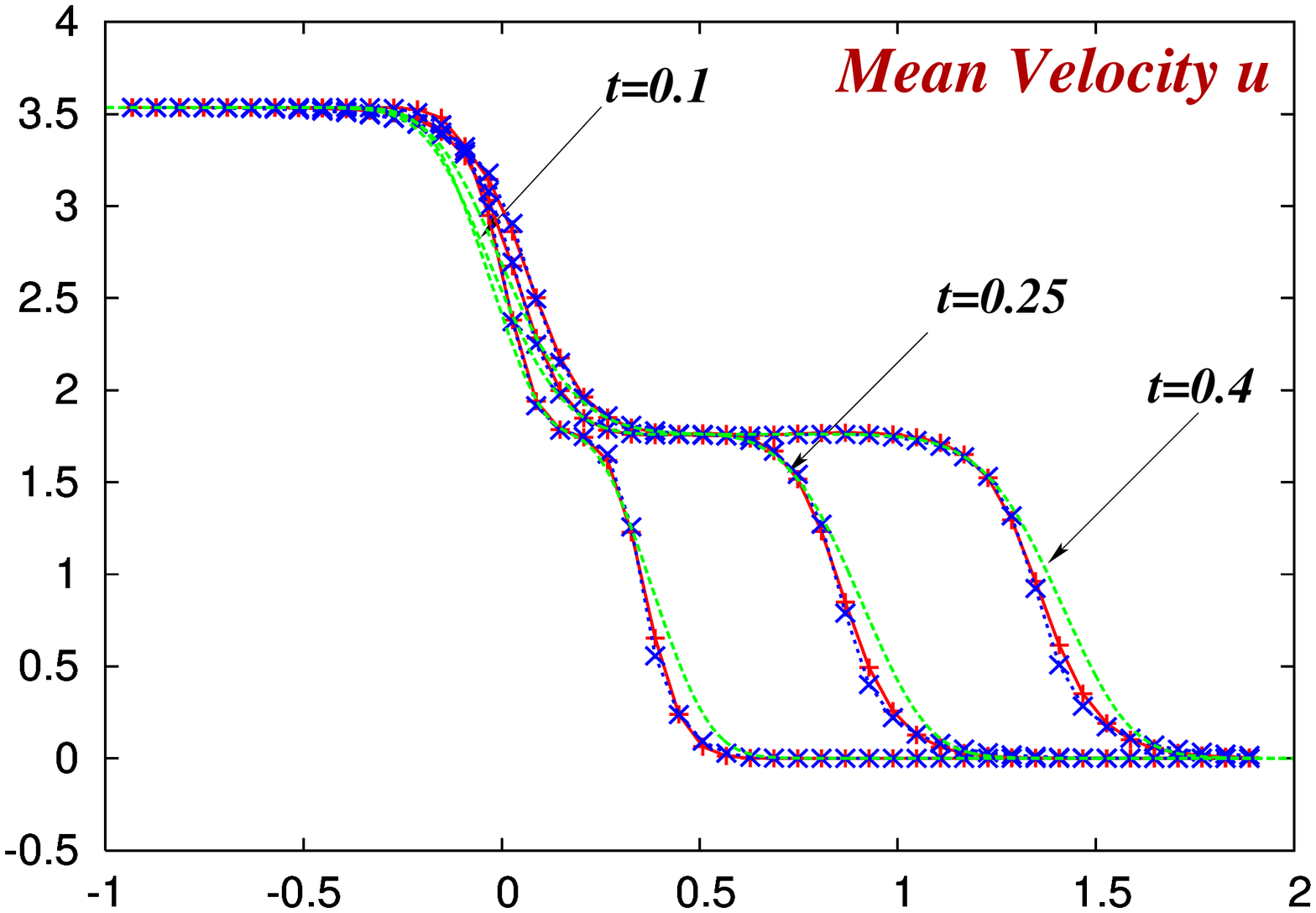}    
\\
(1)&(2)
\\
\includegraphics[width=7.5cm]{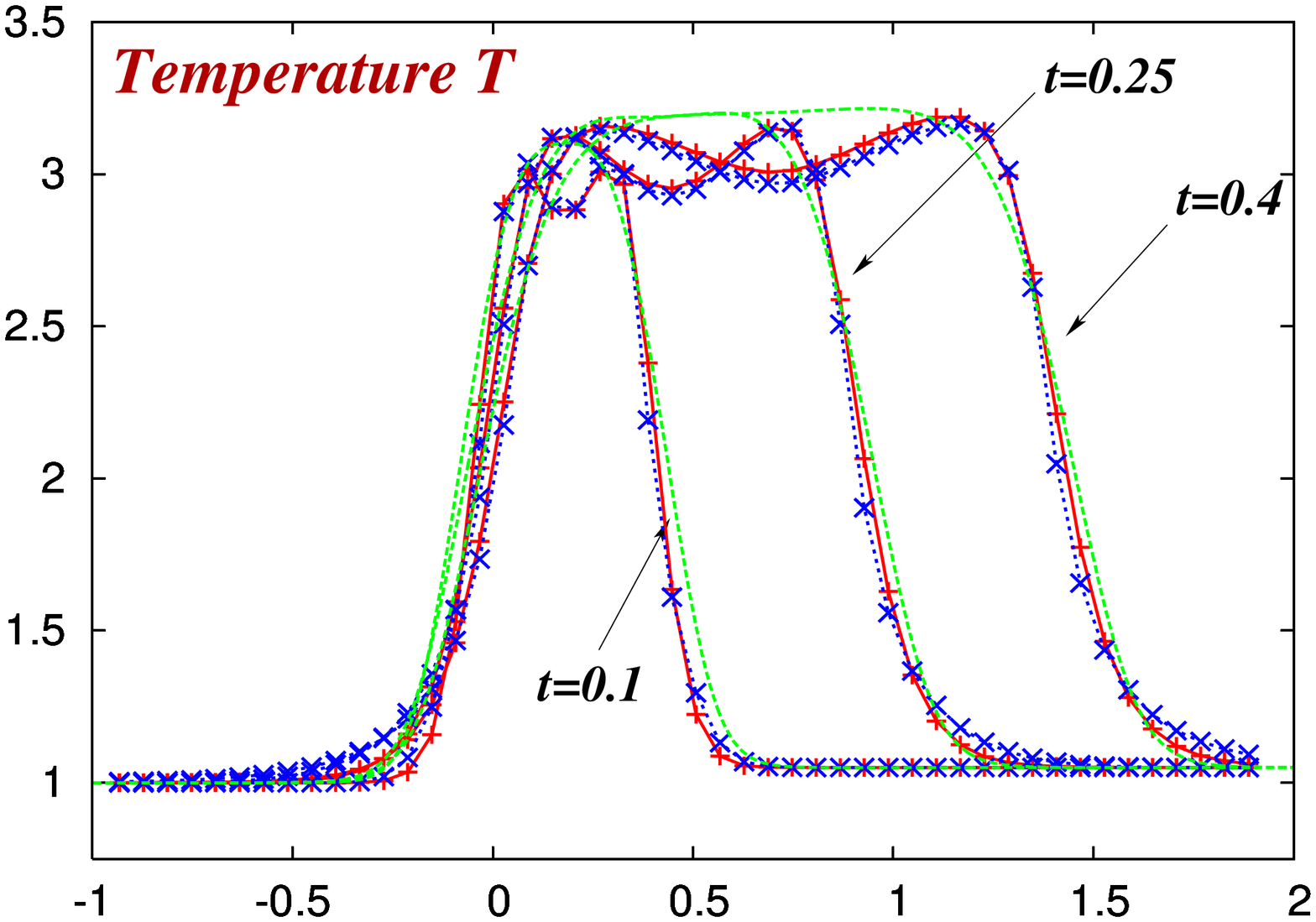}   
&
\includegraphics[width=7.5cm]{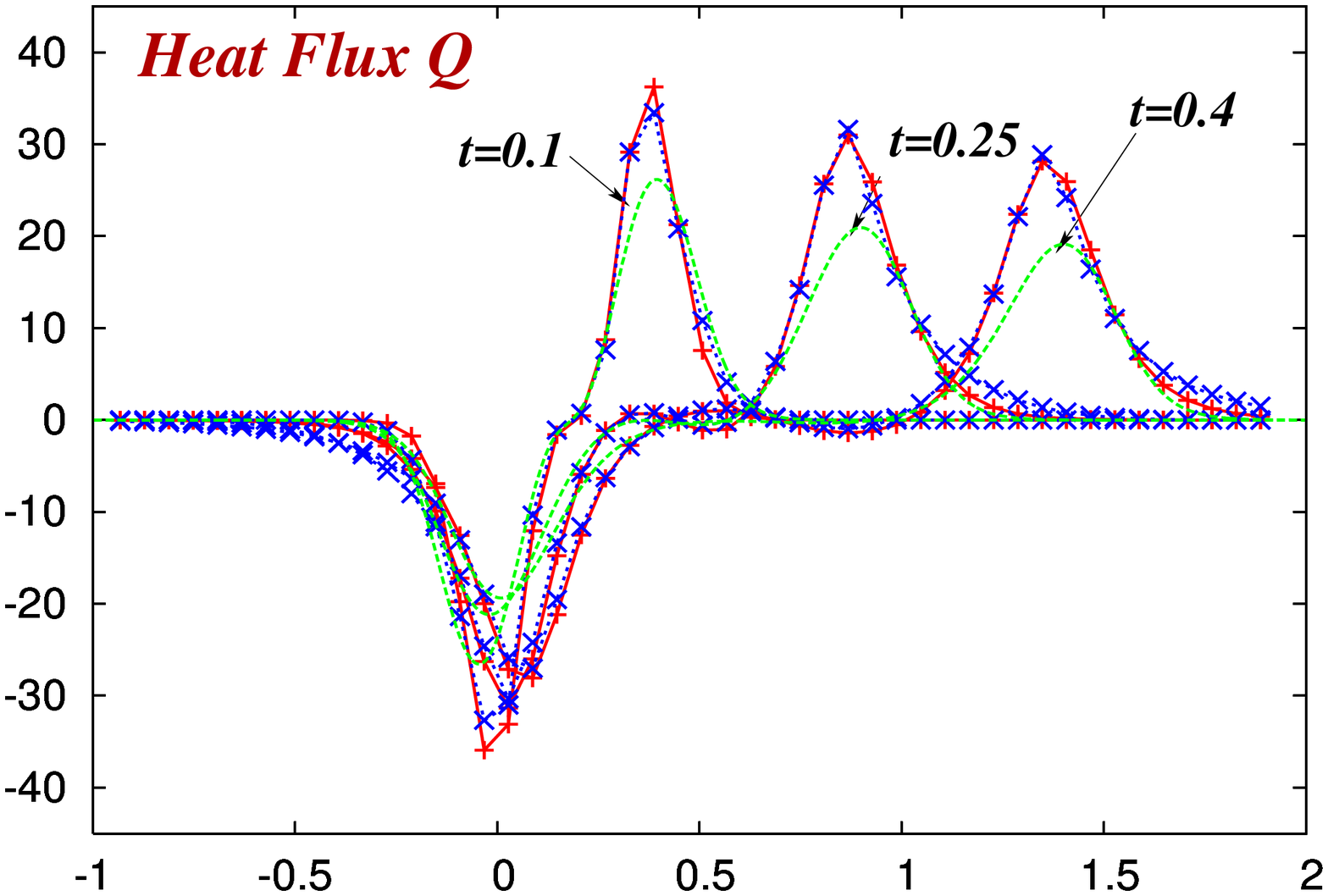}   
\\
(3)&(4)
\end{tabular}
\caption{Riemann problem ($\varepsilon=10^{-1}$): crosses ({\tt +}) represent the numerical solution to the Boltzmann equation  obtained with our method (\ref{AP-1}), stars ({\tt x}) represent the solution  to the ES-BGK model  and lines is the solution corresponding to the compressible Navier-Stokes system.  Evolution of  (1) the density $\rho$, (2) mean velocity $u$, (3) temperature $T$ and  (4) heat flux $\QQ$ at time $t=0.1$, $0.25$ and $0.4$.}
\label{fig:04-1}
\end{figure}

 Therefore, in  Figure~\ref{fig:04-1}, we report the numerical results for  $\varepsilon=10^{-1}$ and make comparison between the numerical solution obtained with the scheme (\ref{AP-1}) and the one obtained with a high order explicit method for the compressible Navier-Stokes. In this case, the behavior of macroscopic quantities (density, mean velocity, temperature and heat flux) agree very well even if the time step is at least ten times larger with our method (\ref{AP-1}).

\begin{figure}[htbp]
\begin{tabular}{cc}
\includegraphics[width=7.5cm]{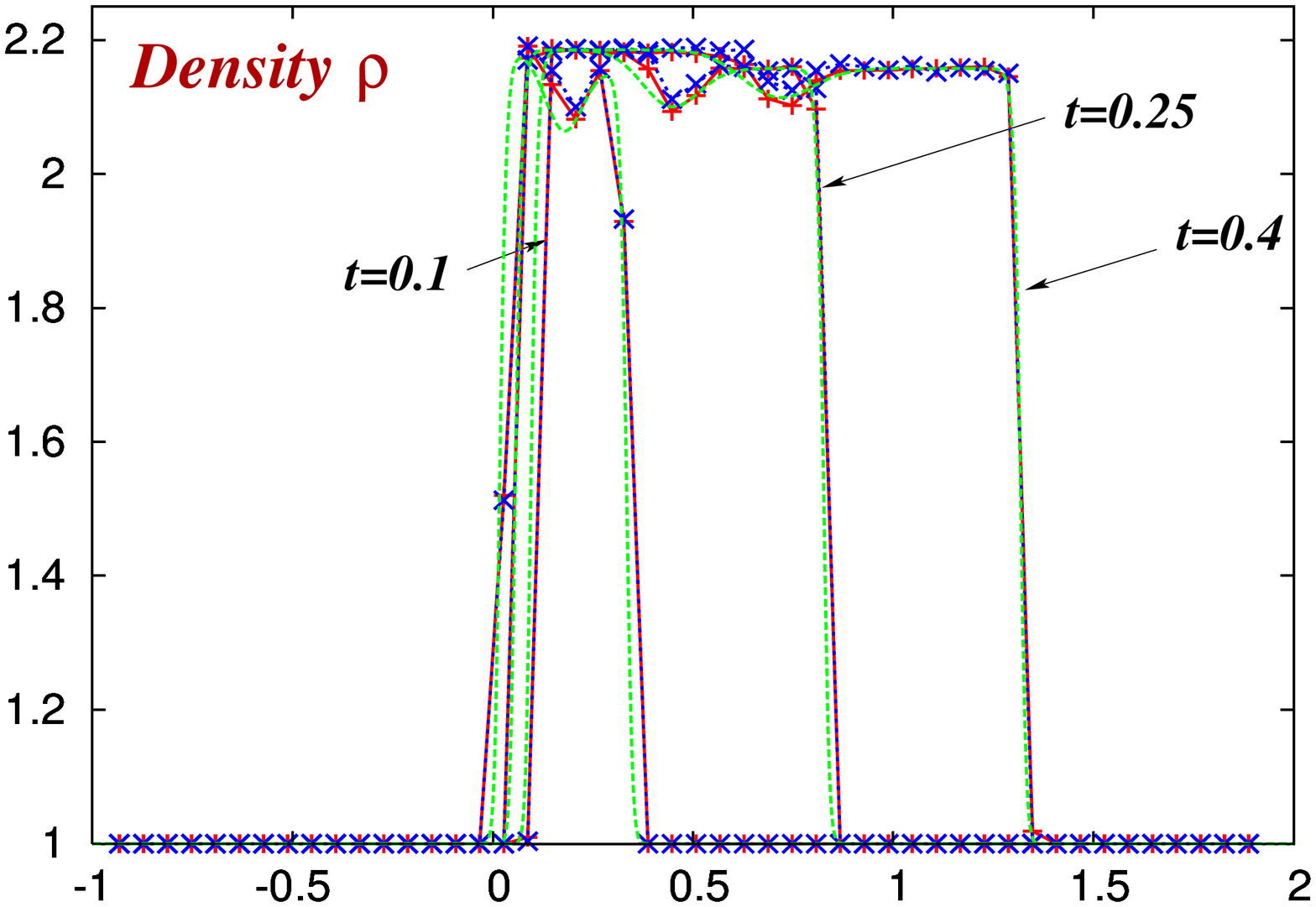}    
&
\includegraphics[width=7.5cm]{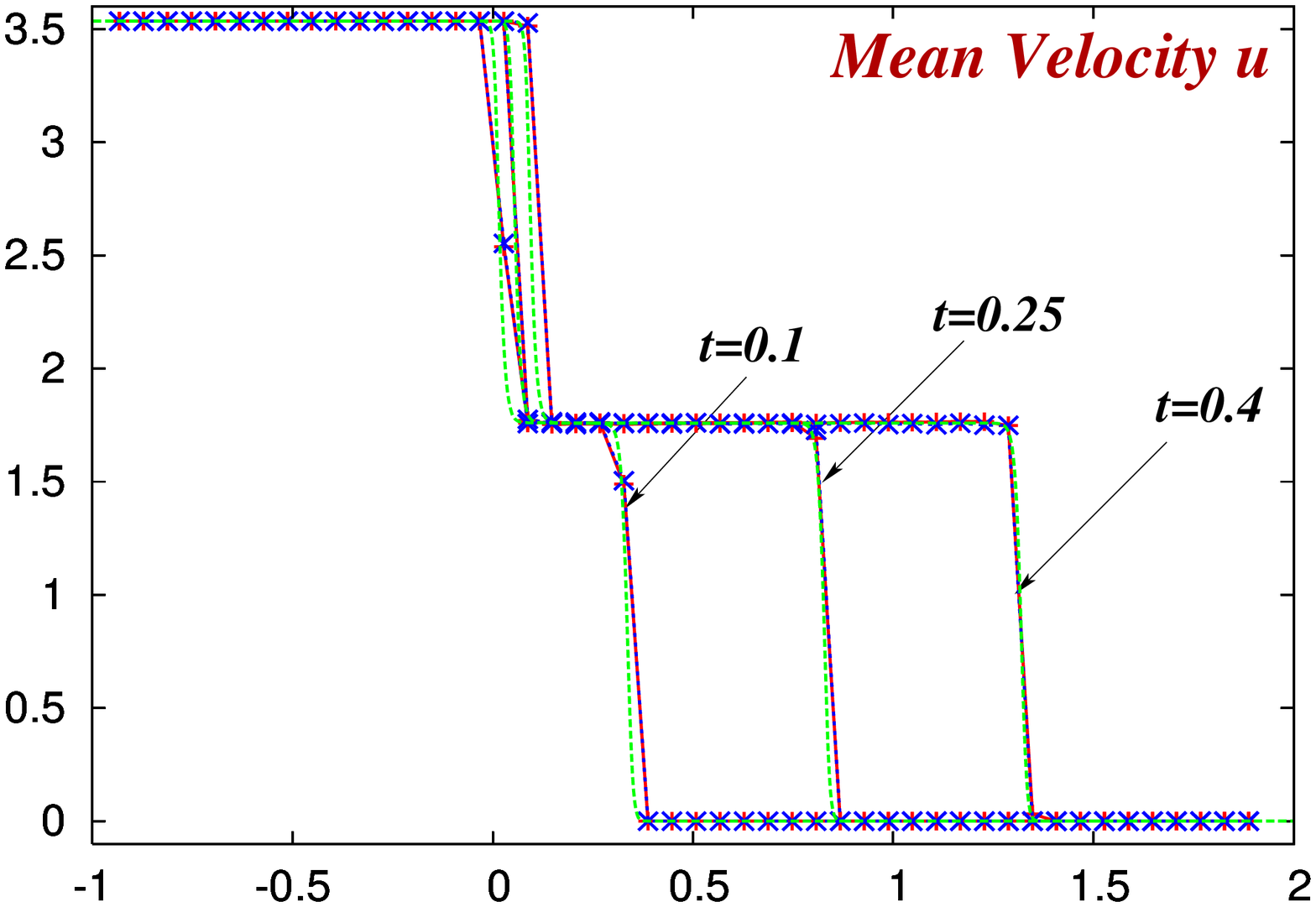}    
\\
(1)&(2)
\\
\includegraphics[width=7.5cm]{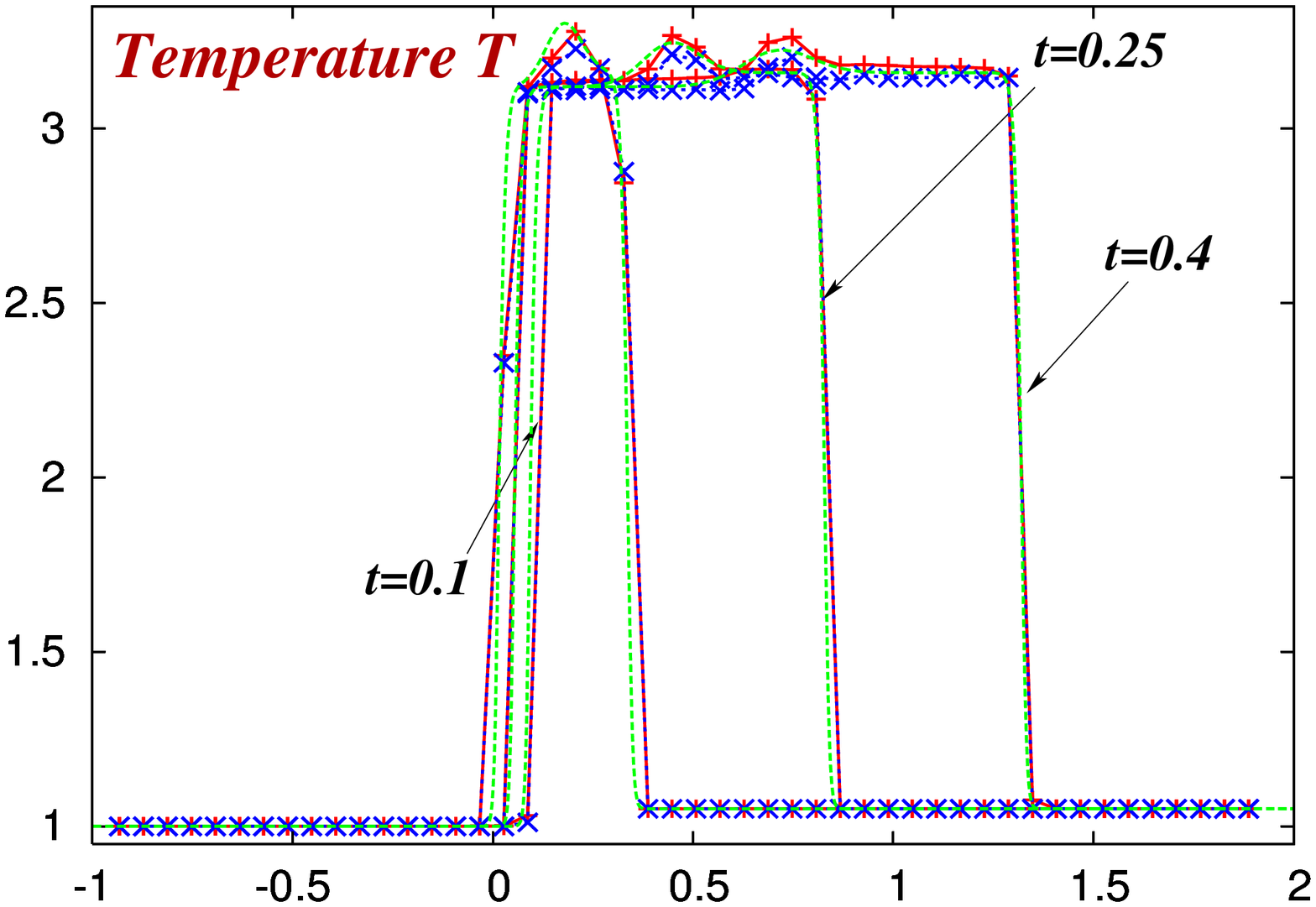}   
&
\includegraphics[width=7.5cm]{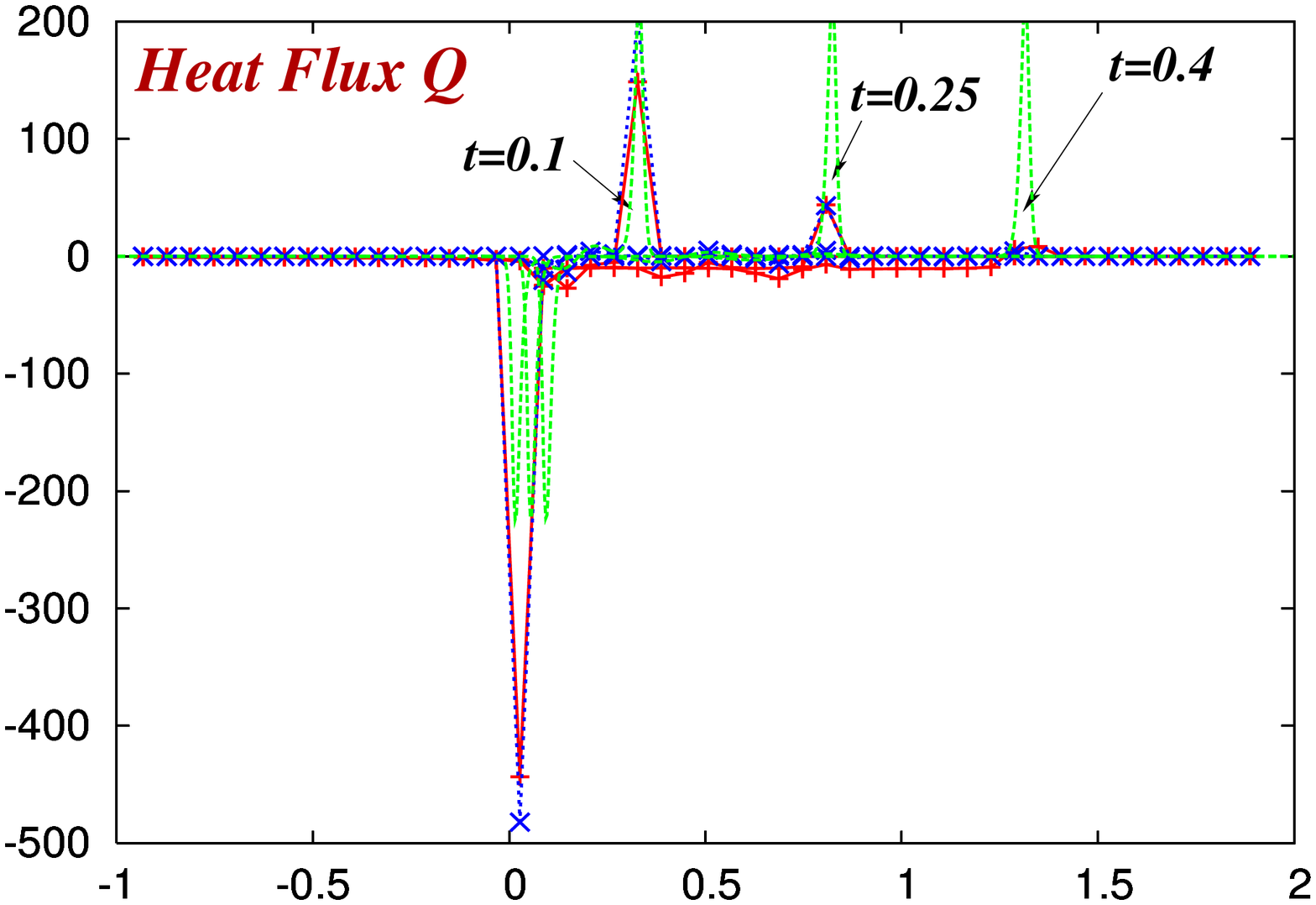}   
\\
(3)&(4)
\end{tabular}
\caption{Riemann problem ($\varepsilon=10^{-3}$): crosses ({\tt +}) represent the numerical solution to the Boltzmann equation  obtained with our method (\ref{AP-1}), stars ({\tt x}) represent the solution  to the ES-BGK model  and lines is the solution corresponding to the compressible Navier-Stokes system. Evolution of  (1) the density $\rho$, (2) mean velocity $u$, (3) temperature $T$ and  (4) heat flux $\QQ$ at time $t=0.1$, $0.2$ and $0.3$.}
\label{fig:04-2}
\end{figure}

Finally in  Figure~\ref{fig:04-2},  we choose $\varepsilon=1.\,10^{-3}$. In this problem, the density, mean velocity and temperature are relatively close to the one obtained with the approximation of the Navier-Stokes system. Even the qualitative behavior of the heat flux  agrees well with  the heat flux  corresponding to the compressible Navier-Stokes  system $\kappa\,\nabla_x T$, with $\kappa=\rho\,T$ (see Figure~\ref{fig:04-2}), yet some differences can be observed, which means that the use of ES-BGK models to derive macroscopic models has a strong influence on the heat flux.

\subsection{Flow around a cylinder}
This example has been considered in for instance  \cite{bib:21}. The computational domain is set to be $[-20, 20] \times [-20, 20]$. The cylinder is centered at the origin, with a diameter of $2$.

We consider an incoming flow at the boundary $\|x\|=8$  with the following conditions: $\rho_{i} = 1$, $u_i = (M\sqrt{2T_i}, 0)^T$ , $T_i = 1$ with $M=0.1$ and $M=0.5$. The freestream Knudsen number ranges from $\varepsilon = 0.1$ to $\varepsilon=10^{-3}$. 

Concerning boundary condition, the wall of cylinder is considered as a pure diffusive boundary conditions
$$
f(t,x,v) = \frac{\rho(t,x)}{(2\pi\,T_w)^{d_v/2}} \,\exp\left(-\frac{|v|^2}{2\,T_w}\right), \quad v\cdot n_x \,<\, 0,\, \quad{\rm and}\quad \|x\| = 1,
$$
where $\rho$ is computed such that the global flux is zero at the boundary (and mass is preserved) and $T_w=1.05$. 

 To start the calculation take a uniform initial solution equal to the values defined by the boundary conditions:
$$
f_0(x,v) = \frac{\rho_i}{(2\pi\,T_i)^{d_v/2}} \,\exp\left(-\frac{|v-u_i|^2}{2\,T_i}\right), \quad v\in \R^2, \quad 1 \leq \|x\|\leq 8.
$$
 Then, we solve the kinetic  equations for the different grid densities considered, until a steady state is reached.

We define the Mach number from the macroscopic quantities, computing the moments of the  distribution function with respect to $v\in\R^2$, by 
$$
M^2 = \frac{\|u\|^2}{\gamma\, T}, 
$$
where $c:=\sqrt{\gamma\,T}$ is the sound speed. 

We apply our numerical scheme (\ref{AP-1}) to the ES-BGK equation and plot the numerical results in the following figures (Fig. \ref{sphere:0},\ref{sphere:1},\ref{sphere:2})  and the solution can be compared to the numerical solution of the Euler equations \cite{bib:21}.

Figure \ref{sphere:0} shows the density contours at different times for a free streaming Mach number $M=0.1$ whereas the Knudsen number is $\varepsilon=0.01$. The reflecting shock, and the Mach shock can all be identified. However, the reflecting shock is not identified as kinetic, since both density and temperature are large on the shock, which makes its distribution much closer to the Maxwellian than the other two shocks. Also, the two seperation points, where the gas is the most rarefied, are well captured and correctly identified as kinetic. In Figure \ref{sphere:1}, the contour plot of the local Mach number, the density and the temperature is shown when the steady state is reached.

\begin{figure}[htbp]
\begin{tabular}{cc}
\includegraphics[width=7.5cm]{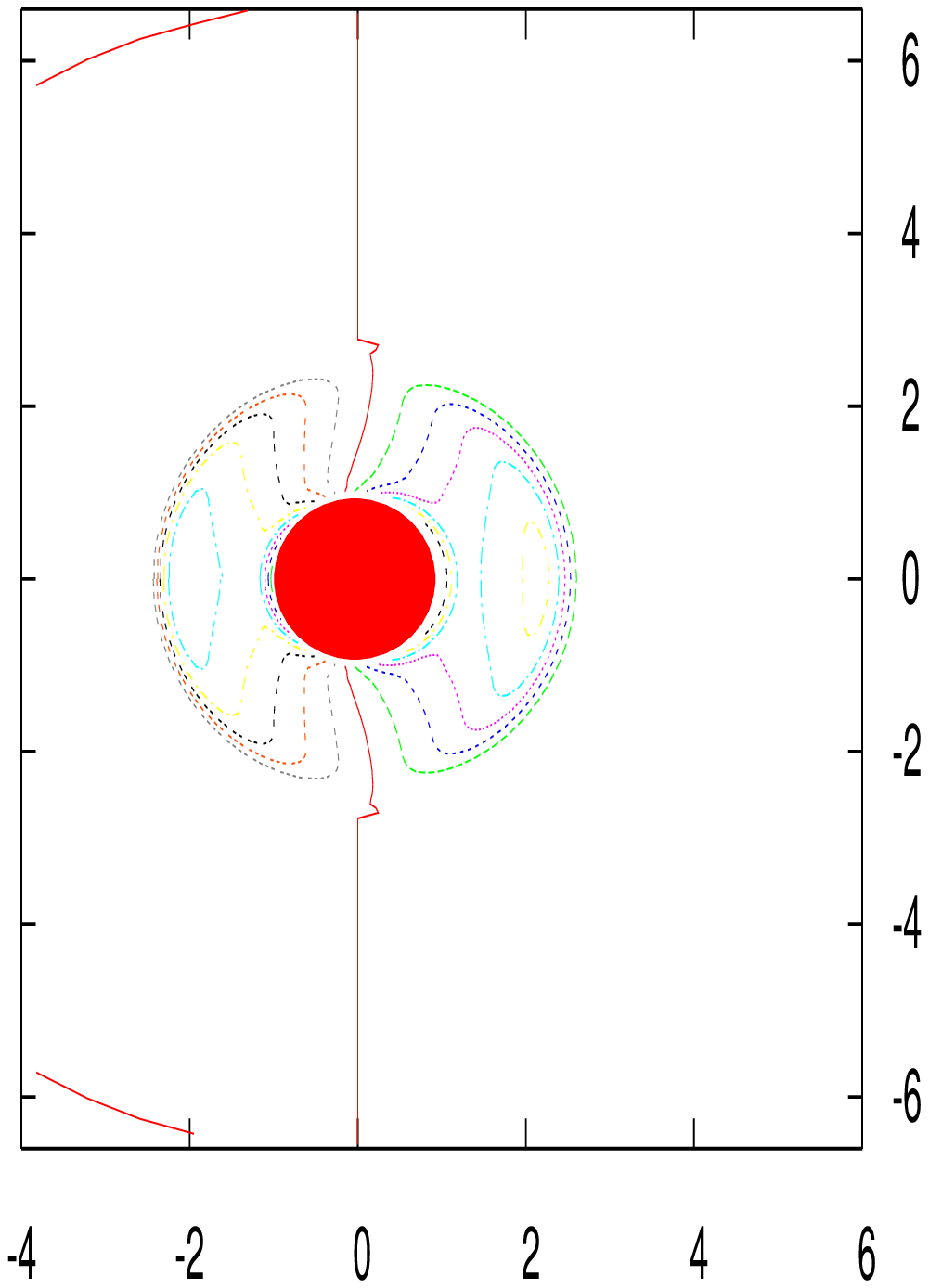}    
&
\includegraphics[width=7.5cm]{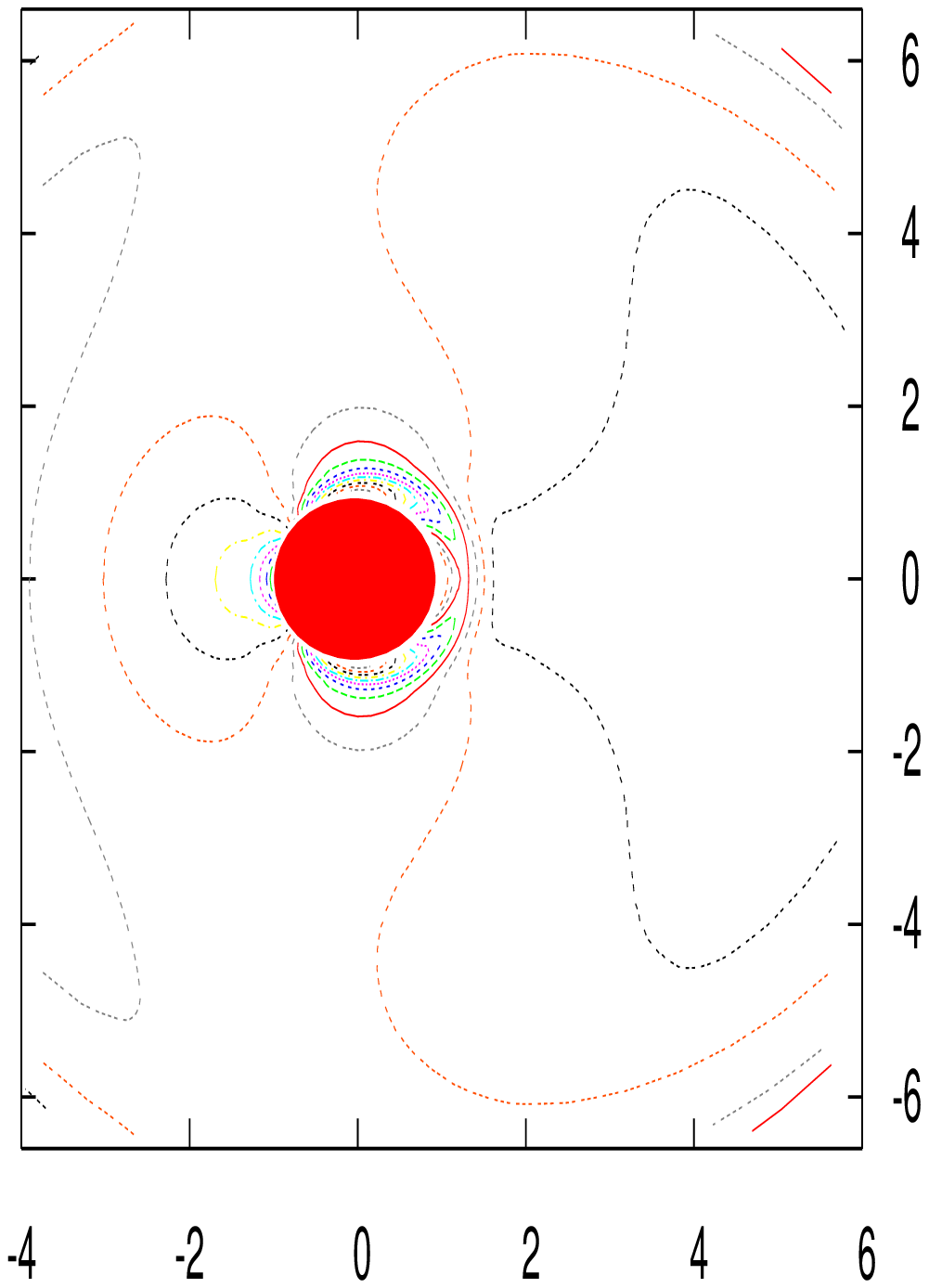}    
\\
(1)&(2)
\\
\includegraphics[width=7.5cm]{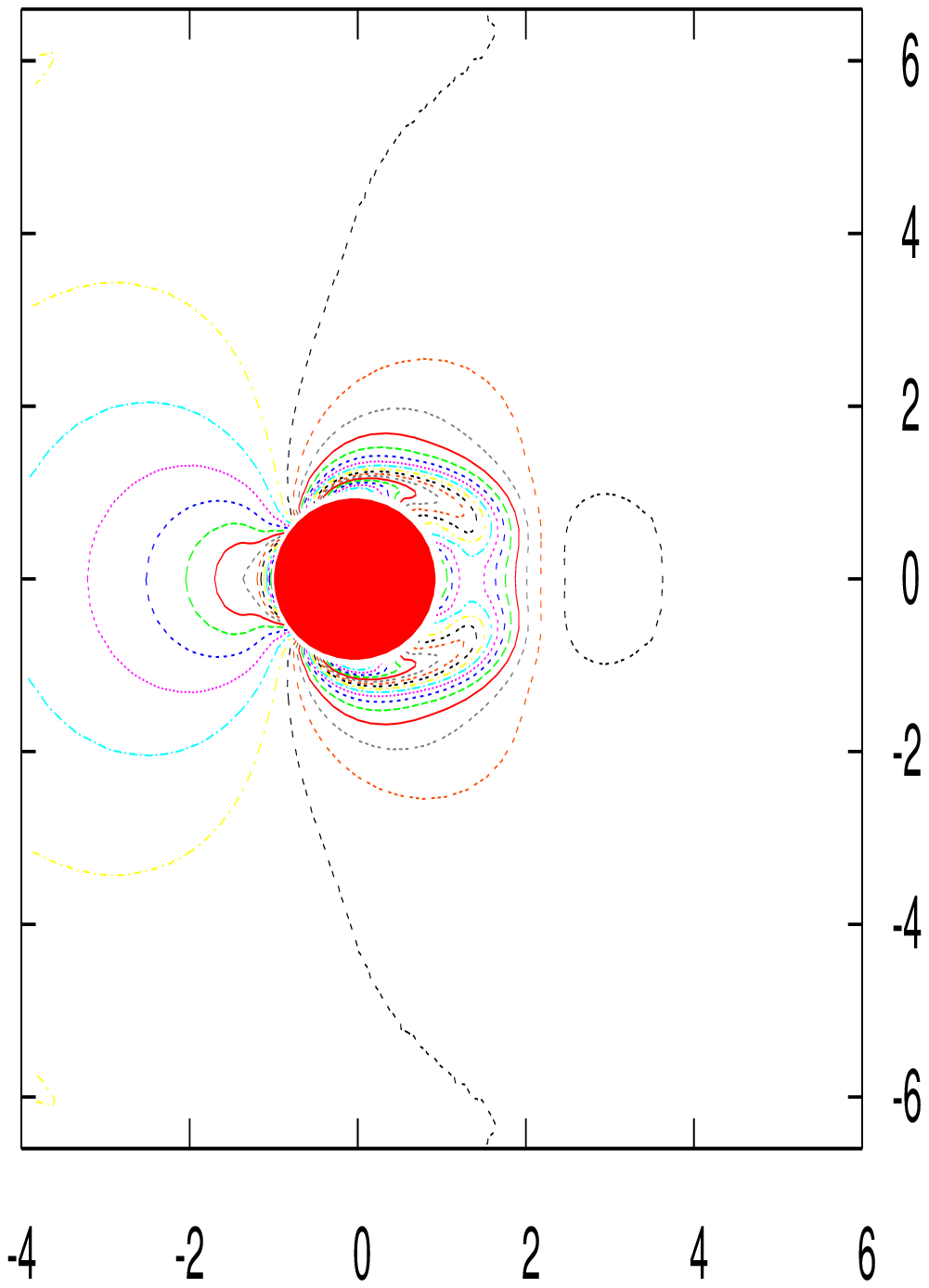}
&
\includegraphics[width=7.5cm]{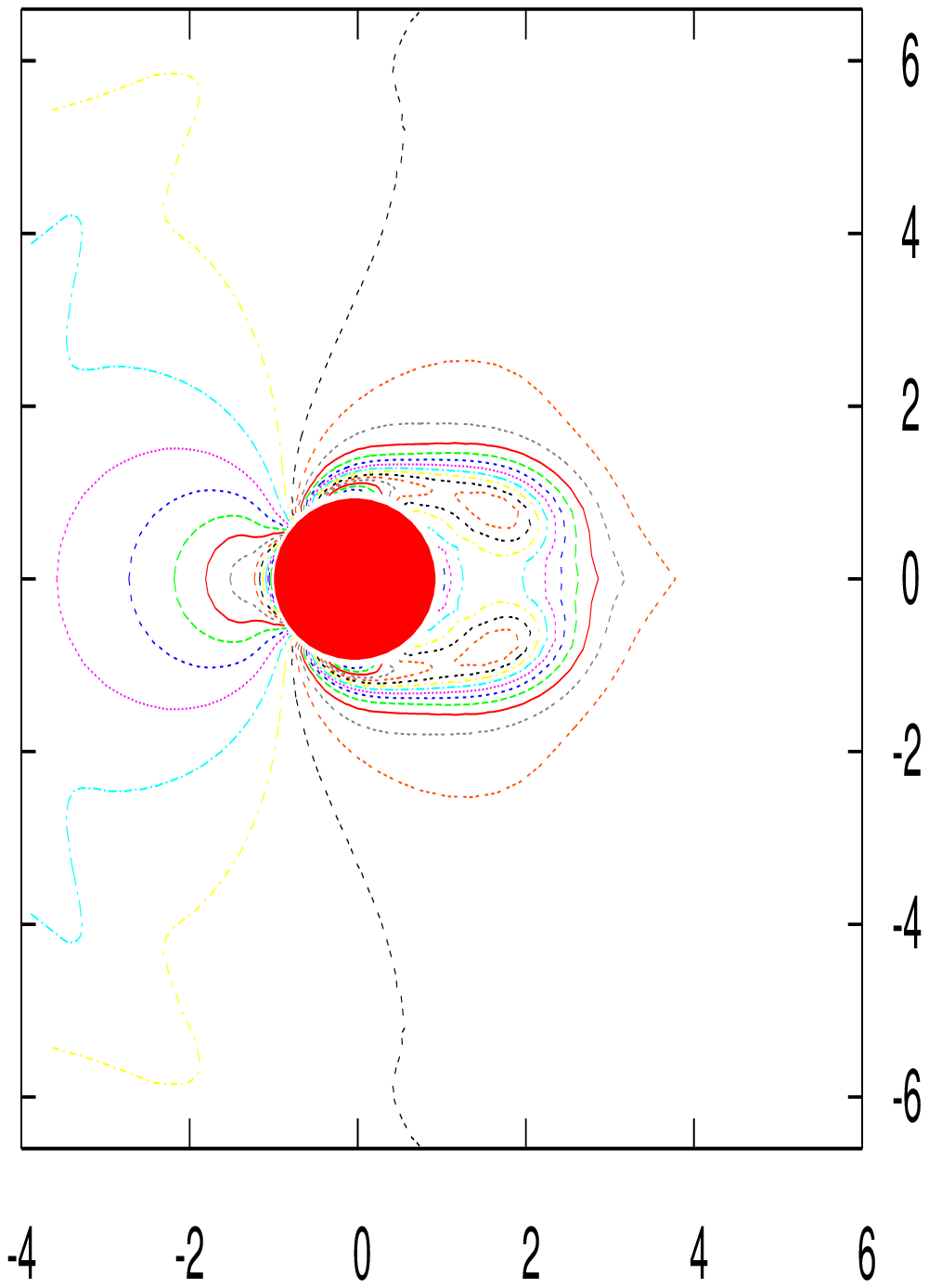}
\,   
\\
(3) &\,
\end{tabular}
\caption{Flow around a cylinder ($M=0.1$, $\varepsilon=0.01$), the numerical solution of the ES-BGK model obtained with our method (\ref{AP-1}) : evolution of the density contour  at time (1) $t=1$  (2) $t=6$ (3)  $t=16$ (4)  $t=30$.
}
\label{sphere:0}
\end{figure}

\begin{figure}[htbp]
\begin{tabular}{ccc}
\includegraphics[width=5.125cm]{test_sphereMach0.1epsi_0.01rhosteady.eps}    
&
\includegraphics[width=5.125cm]{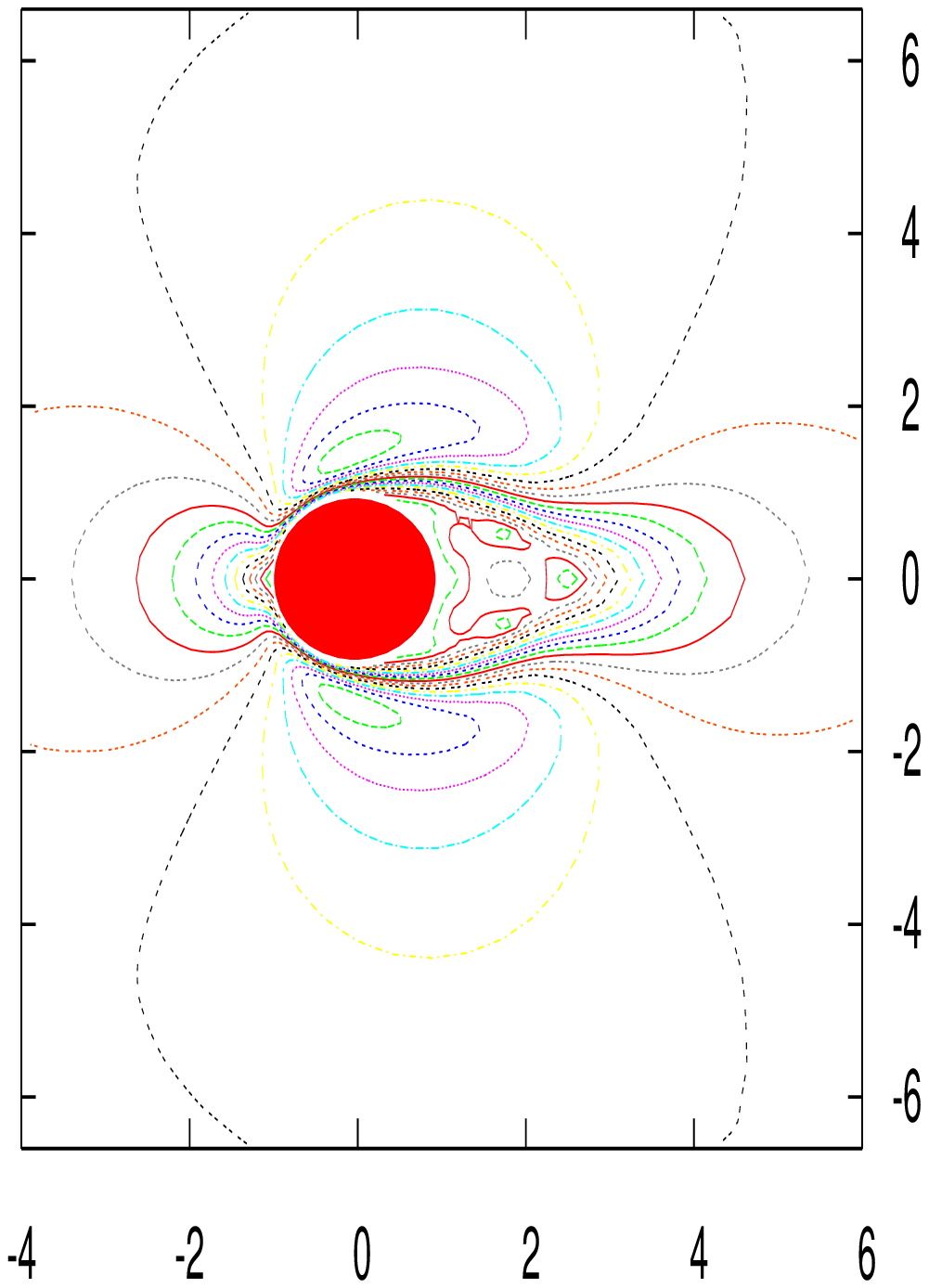}    
&
\includegraphics[width=5.125cm]{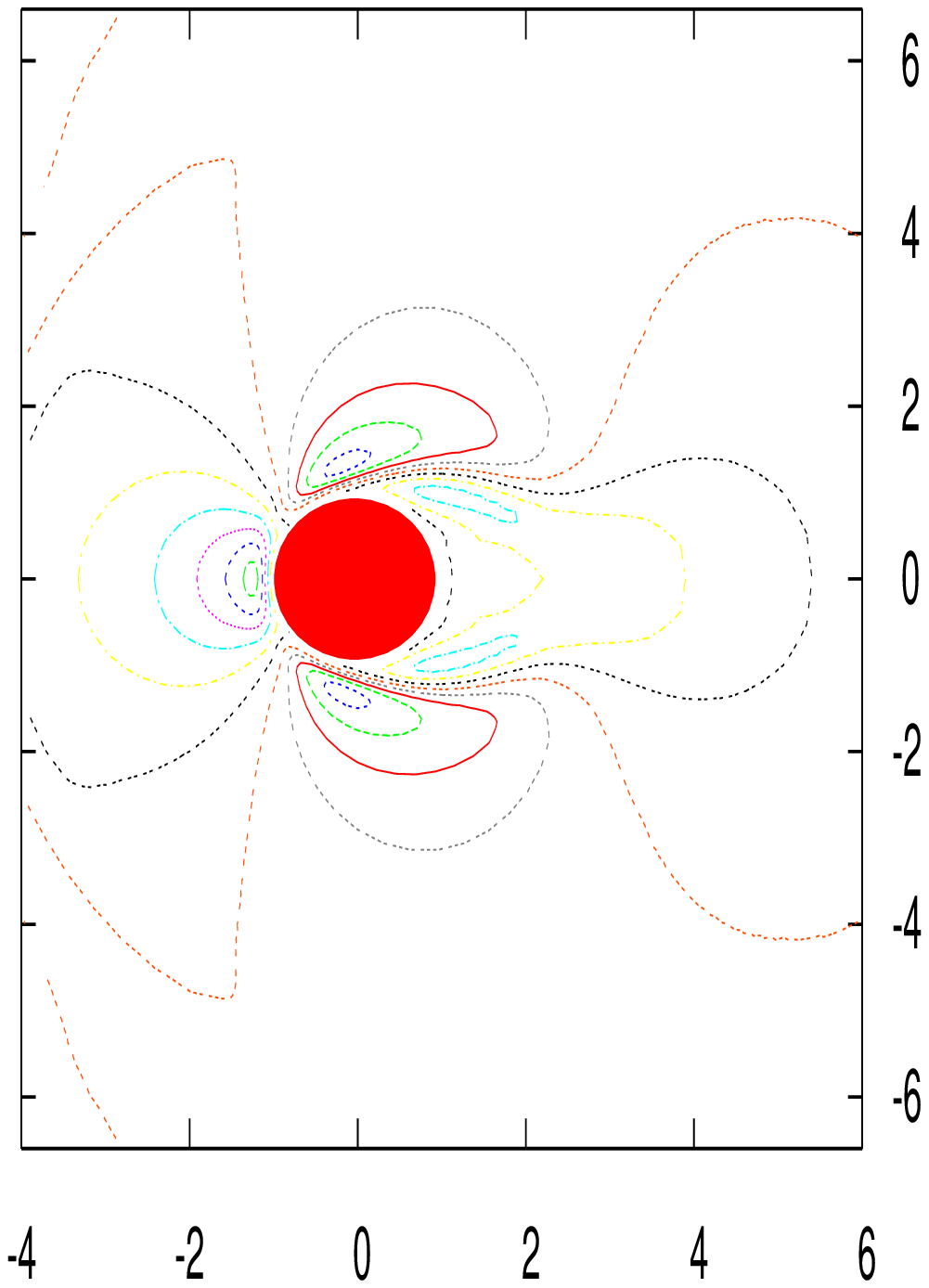}    
\\
(1)&(2)&(3)
\end{tabular}
\caption{Flow around a cylinder ($M=0.1$, $\varepsilon=0.01$), the numerical solution of the ES-BGK model obtained with our method (\ref{AP-1}) :steady state (1) density $\rho$  (2) local Mach number $M$ (3)  local temperature $T$.
}
\label{sphere:1}
\end{figure}

\begin{figure}[htbp]
\begin{tabular}{ccc}
\includegraphics[width=5.125cm]{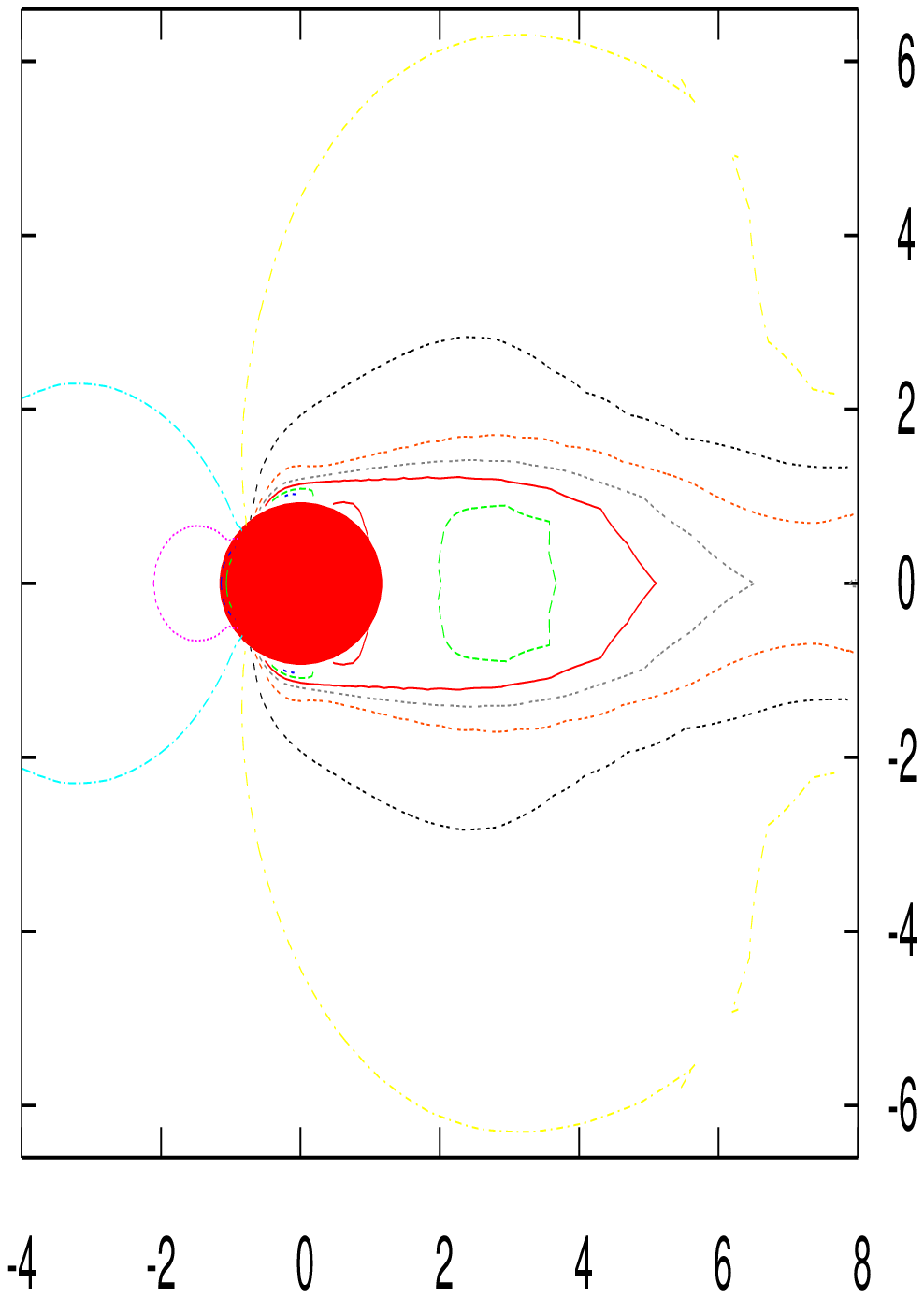}    
&
\includegraphics[width=5.125cm]{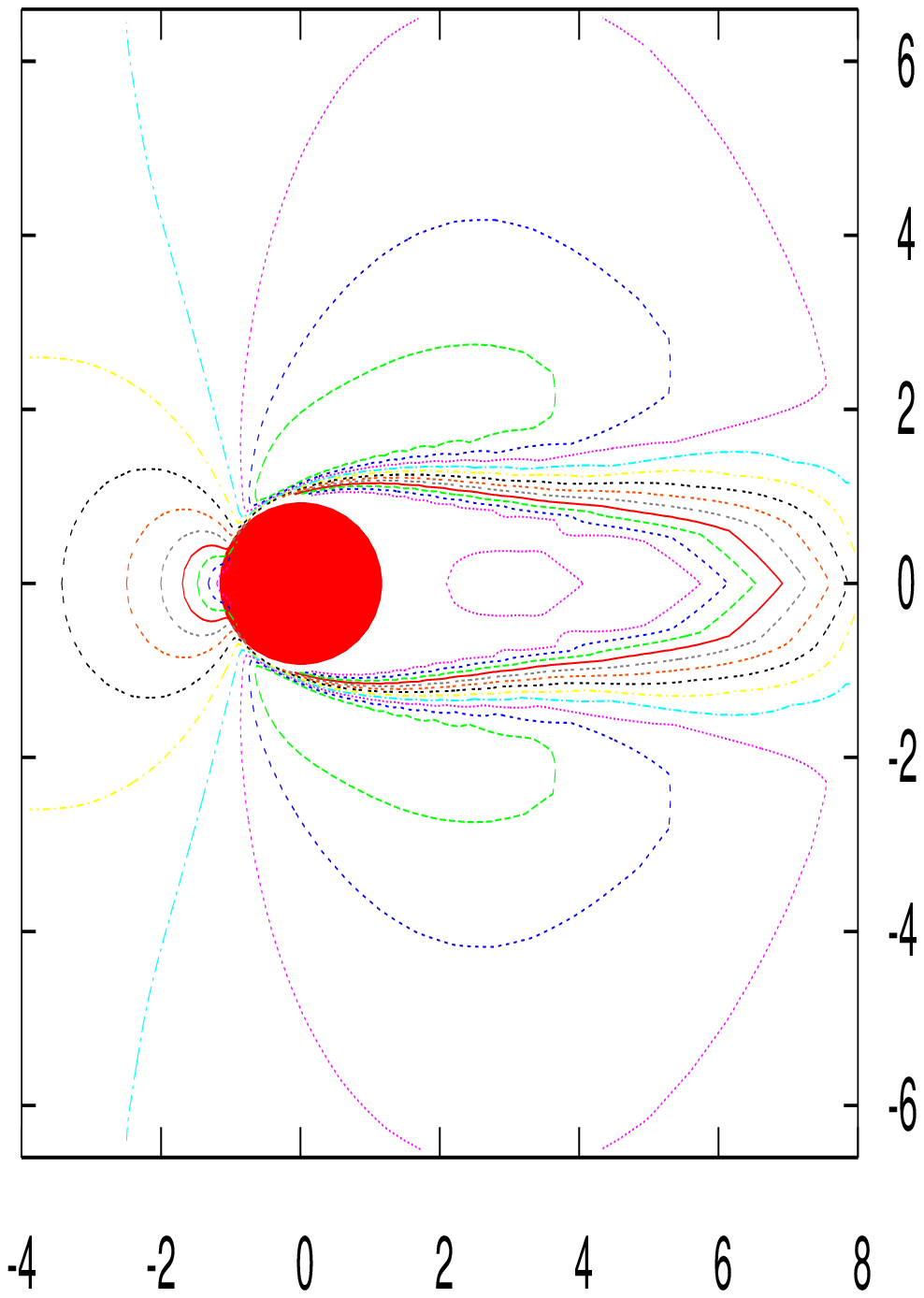}    
&
\includegraphics[width=5.125cm]{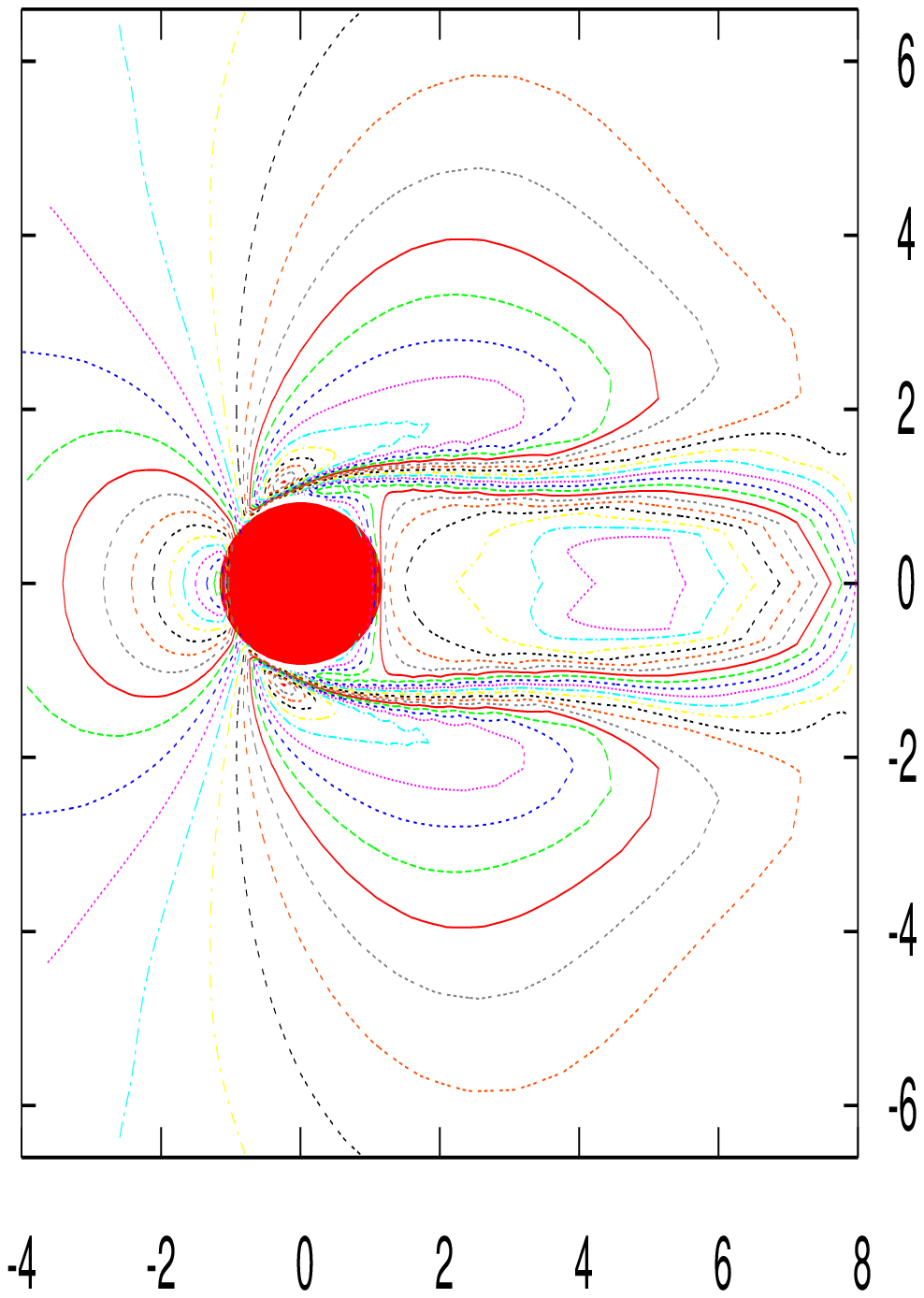}    
\\
(1)&(2)&(3)
\end{tabular}
\caption{Flow around a cylinder ($M=0.5$, $\varepsilon=0.01$), the numerical solution of the ES-BGK model obtained with our method (\ref{AP-1}) :steady state (1) density $\rho$  (2) local Mach number $M$ (3)  local temperature $T$.
}
\label{sphere:2}
\end{figure}

\begin{figure}[htbp]
\begin{tabular}{cc}
\includegraphics[width=7.5cm]{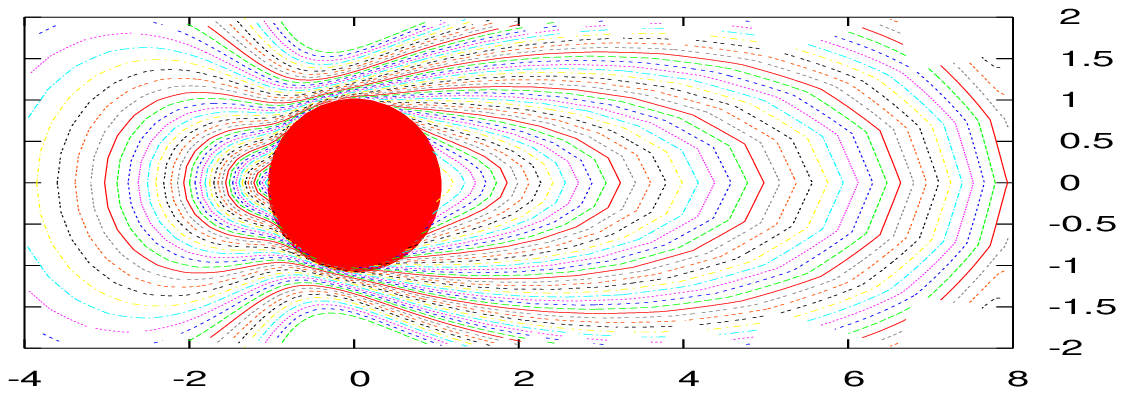}    
&
\includegraphics[width=7.5cm]{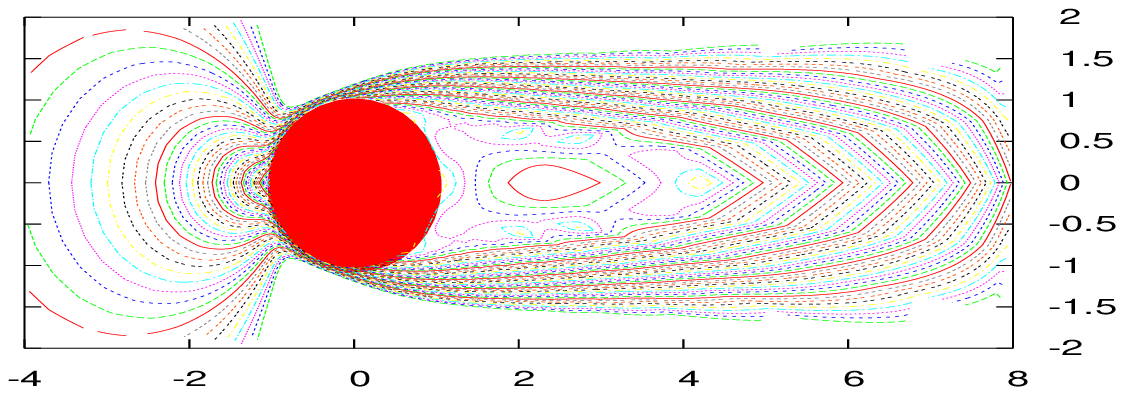}    
\\
(1)&(2)
\\
\includegraphics[width=7.5cm]{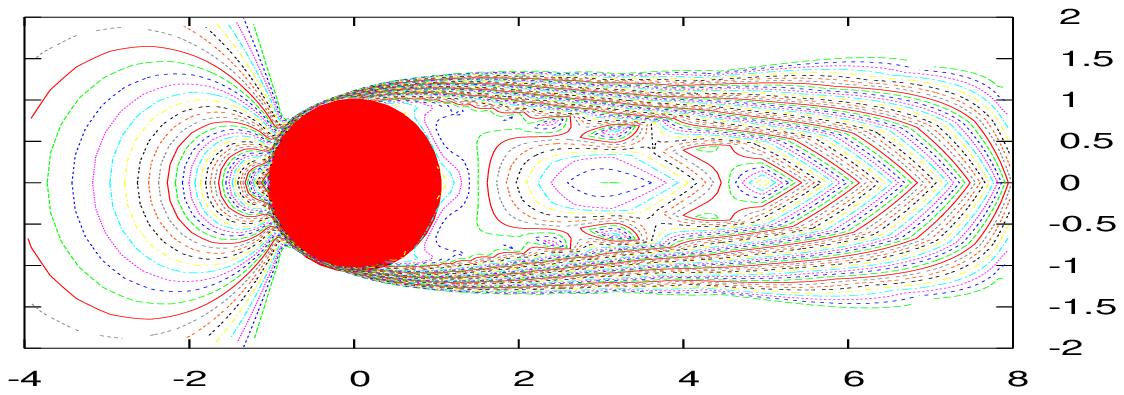}    
&
\includegraphics[width=7.5cm]{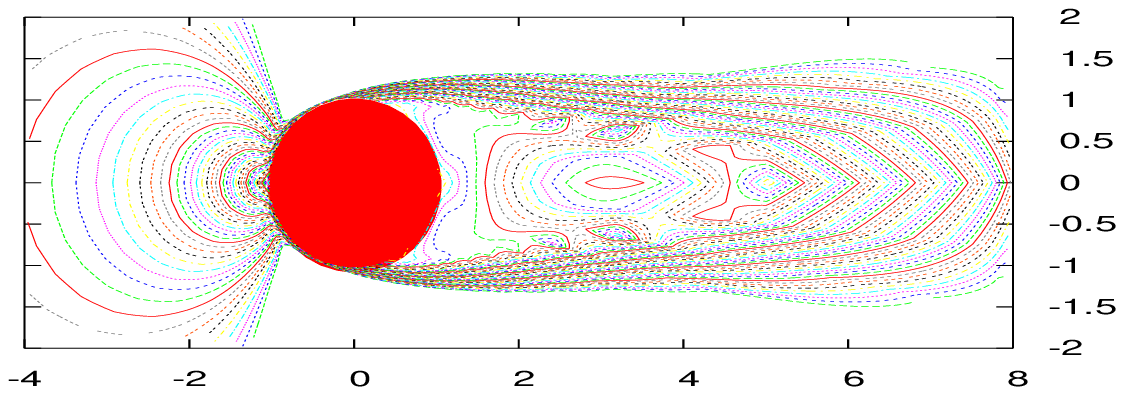}    
\\
(3)&(4)
\end{tabular}
\caption{Flow around a cylinder ($M=0.5$), the numerical solution of the ES-BGK model obtained with our method (\ref{AP-1}) : steady state of the local Mach number  for different values of Knudsen number (1) $\varepsilon=0.5$  (2) $\varepsilon=0.1$, (3)  $\varepsilon=0.01$,  (4)  $\varepsilon=0.001$.
}
\label{sphere:3}
\end{figure}

\section{Conclusion}

In this paper we present an accurate deterministic method for the numerical
approximation of the space inhomogeneous, time dependent ES-BGK equation. The method is a temporally implicit-explicit scheme to deal with the stiffness of the collision operator. The computational cost of the implicit part is close
to an explicit one, without using any nonlinear algebraic system solver,
by utilizing the particular structure of the ES-BGK operator. This effective time discretization allows the treatment of problems with a broad
 range of mean free path. Moreover, the numerical results, and the comparison with other techniques, show the effectiveness of the present method for a wide class of problems.


\begin{flushleft} 
\signff 
\end{flushleft}
\vspace{-4.25cm}
\begin{flushright} 
\signsj 
\end{flushright}

\end{document}